\documentclass[12pt]{article}

\usepackage{graphicx}
\usepackage{pst-plot}
\usepackage{pstricks-add}
\usepackage{amssymb}
\usepackage{url}
\usepackage[letterpaper,margin=1in]{geometry}

\begin{document} 

\newtheorem{theorem}{Theorem}[section]
\newtheorem{korollar}{Corollary}[section]
\newtheorem{remark}{Bemerkung}[section]
\newtheorem{aufgabe}{Aufgabe}[section]
\newtheorem{define}{Definition}[section]
\newtheorem{example}{Example}[section]
\newtheorem{lemma}{Lemma}[section]
\begin{center}
{\Large \bf On the Shape of the Symmetric Solution Set of a Linear Complementarity Problem with Interval Data}\\[2ex]
{\bf Uwe Sch\"afer}\\[1ex]
Josef-Durler-Schule, Technisches Gymnasium\\
Richard-Wagner-Ring 24\\
76437 Rastatt\\
E-Mail: Dr.Schaefer@jdsr.de
\end{center}

\begin{center}
{\bf Abstract}
\end{center}
In this paper we give some two-dimensional and some three-dimensional examples for the shape of the symmetric solution set of a linear complementarity problem where the given data are not explicitly known but can only be enclosed in intervals.

\tableofcontents

\section{Introduction}

Given an $n \times n$ matrix $M$ and an $n$-vector $q$ the linear complementarity problem, abbreviated LCP,  involves finding two vectors $w$ and $z$ such that
\begin{equation} \label{defa}
w-Mz=q, \quad w\geq o, \, z\geq o, \quad w^{\mbox{\scriptsize T}}z=0
\end{equation}
or showing that no such vectors exist. Here, matrix and vector inequalities are meant componentwise and $o$ denotes the zero vector. Due to the big number of applications there is also another definition of the LCP which is equivalent to (\ref{defa}). Given an $n \times n$ matrix $M$ and an $n$-vector $q$ the linear complementarity problem  involves finding a vector $z$ such that
\begin{equation} \label{defb}
q+Mz\geq o, \quad  z\geq o, \quad (q+Mz)^{\mbox{\scriptsize T}}z=0
\end{equation}
or showing that no such vector exists. For a detailed introduction to the LCP we refer to the books \cite{cottlepangstone}, \cite{murty}, and \cite{schaeferbuch}.

In \cite{schaeferdiss} we have already considered the problem that the matrix $M$ and the vector $q$ are not explicitly known but can only be enclosed in intervals. Such a problem arises, for example, from discretizing an ordinary free boundary-value problem without neglecting the remainder term of the underlying Taylor expansion, see \cite{schaeferNFAO} and \cite{schaeferdiss}. 

So, given an interval matrix $[M]$ and an interval vector $[q]$ we are interested in the solution set
\[
\Sigma([M],[q])=
\left\{
z  \in \mathbb{R}^n:  \mbox{ there exist $M \in [M]$ and $q \in [q]$ satisfying (\ref{defb}) } 
\right\}.
\]
In \cite{aleschaefer}, \cite{donauer}, \cite{schaeferhabil}, and \cite{schaeferdiss} the solution set was already studied and visualized for the case $n=2$ and $n=3$. In this paper we consider the symmetric solution set
 \[
\Sigma_{sym}([M],[q])=
\left\{
z  \in \mathbb{R}^n:  \mbox{ there exist $M \in [M]$ and $q \in [q]$ satisfying (\ref{defb}) and $M=M^{\mbox{\scriptsize T}}$ } 
\right\}.
\]
The paper is organized as follows. After the preliminaries we present a new theorem that 
generalizes Theorem 3.2 in \cite{aleschaefer}. Then we present, in general, the way how to describe the solution set  $\Sigma([M],[q])$ and the 
symmetric solution set $\Sigma_{sym}([M],[q])$. Concerning the latter we use the Fourier-Motzkin elimination process presented in \cite{alefeldkreinmayer}. Finally, we give several 
examples.

\section{Preliminaries}
\subsection{Interval arithmetic}
We consider compact intervals
\[
[\underline{a},\overline{a}]:=\{ x \in \mathbb{R}: \underline{a} \leq x \leq \overline{a} \} 
\]
and denote the set of all such intervals by ${\bf IR}$. We also write $[a]$ instead of $[\underline{a},\overline{a}]$. Furthermore, we consider matrices with a compact interval in each of its elements; i.e., $[\underline{A},\overline{A}]=([a_{ij}])=([\underline{a}_{ij},\overline{a}_{ij}])$ and
\[
[\underline{A},\overline{A}]=\{ A \in \mathbb{R}^{m \times n}: \underline{A} \leq A \leq \overline{A}\},
\]
and we also write $[A]$ instead of $[\underline{A},\overline{A}]$. By ${\bf IR}^{m \times n}$ we denote the set of all these so-called interval matrices. The set of interval vectors with $n$ components is defined in the same way and is denoted by ${\bf IR}^{n}$.

The intersection of two interval vectors is understood componentwise. If there is at least one component for which the intersection is empty, then the intersection of the two interval vectors is set $\emptyset$, which is the empty set. For completeness, we recall that for $[a],\, [b]  \in {\bf IR}$ we have
\begin{equation} \label{zentralbeeck}
[a] \cap [b] \not= \emptyset \quad \Leftrightarrow \quad \underline{a} \leq \overline{b} \mbox{ and }  \underline{b} \leq \overline{a}.
\end{equation}
For an introduction to interval computations we refer to \cite{alefeld} and \cite{neumaier}. Finally, we want to note that for $A,\, B \in \mathbb{R}^{m \times n}$
\[
A < B \mbox{ means that } a_{ij} < b_{ij} \mbox{ for all } i=1,...,m; \, j =1,...,n.
\]
 
\subsection{Classes of matrices}

A matrix $M \in \mathbb{R}^{n \times n}$ is called a P-matrix if all its principal minors are positive (see \cite{schaeferSIAM}). The relevance of P-matrices to the LCP is established by the following theorem.

\begin{theorem}
A matrix $M \in \mathbb{R}^{n \times n}$ is a P-matrix if and only if the LCP has a unique solution for all $q \in \mathbb{R}^{n}$.
\end{theorem}
For a proof we refer to Theorem 3.3.7 in \cite{cottlepangstone}. By $Z^{n \times n}$ we define the set of matrices with nonpositive off-diagonal entries; i.e., 
\[
M=(m_{ij}) \in Z^{n \times n} \quad
\Leftrightarrow \quad
m_{ij} \leq 0 \,  \mbox{ if } i\not=j.
\]
Furthermore, let $O$ denote the zero matrix.
\begin{define}
We call $M \in  Z^{n \times n}$ an M-matrix, if $M^{-1}$ exists satisfying $M^{-1}\geq O$.
\end{define}
Sometimes M-matrices are also called K-matrices, see Note 3.13.24 in \cite{cottlepangstone}.

\begin{lemma} \label{lemmafan}
Let $A \in  Z^{n \times n}$. Then the following two conditions are equivalent.
\begin{enumerate}
\item $A$ is an M-matrix.
\item There is $u \in  \mathbb{R}^n$ with $u>o$ such that $Au>o$.
\end{enumerate} 
\end{lemma}
For a proof we refer to Lemma 3.4.1 in \cite{schaeferbuch}.
\begin{korollar} \label{korfan}
Let $M, \, A \in  \mathbb{Z}^{n \times n}$ satisfying $M\leq A$. If $M$ is an M-matrix, then $A$ is an M-matrix, too. In particular, we then have $A^{-1}\leq M^{-1}$.
\end{korollar}
The proof can easily be done by Lemma \ref{lemmafan}, where the supplement is an easy consequence of the definition of an M-matrix.
\begin{define}
We call $M \in    \mathbb{R}^{n \times n}$ an $H_+$-matrix, if all its diagonal entries are positive and if its so-called comparison matrix $<M>=(c_{ij})$ defined by
\[
c_{ij}=
\left\{
\begin{array}{lcr}
|m_{ii}| & \mbox{if} & i=j,\\[1ex]
-|m_{ij}| & \mbox{if} & i \not=j,
\end{array}
\right.
\]
is an M-matrix.
\end{define}
Any M-matrix is an $H_+$-matrix, and any $H_+$-matrix is a P-matrix, but not vice versa. See also \cite{schaeferSIAM}.

\section{A new theorem}
\begin{theorem} \label{neuestheorem}
Let $\hat M, \, \tilde M \in  \mathbb{R}^{n \times n}$ satisfying $\hat M \leq \tilde M $. Furthermore, let $\hat q, \, \tilde q \in  \mathbb{R}^{n}$ satisfying $\hat q \leq \tilde q $. Let
$\tilde z$ be a solution of the LCP defined by $\tilde M$ and $\tilde q$ and let $\hat z$ be a solution of the LCP defined by $\hat M$ and $\hat q$. Then,
\[
\mbox{if } \tilde M \, \mbox{ or } \, \hat M \, \mbox{ is an M-matrix, then } \tilde z \leq \hat z.
\]
In addition, let $\tilde z >o$ and let there be $i_0,\, j_0$ such that $\hat m_{i_0j_0} < \tilde m_{i_0j_0}$ or $\hat q_{i_0}< \tilde q_{i_0}$
\[
\mbox{or}
\]
let $\tilde q =-\tilde M \tilde z$ and let there be $i_0$ such that $\hat q_{i_0}< \tilde q_{i_0}$.
Then, 
\[
\mbox{if } \tilde M^{-1} >O \, \mbox{ or } \, \hat M^{-1}>O, \, \mbox{ then } \tilde z < \hat z.
\]
\end{theorem}
{\it Proof: } Case 1: Let $\hat M$ be an M-matrix.\\[1ex]
Then $\hat z $ is the unique solution of  the LCP defined by $\hat M$ and $\hat q$. Let $\Re =\{i: \tilde z_i=0\}$. Then we define $A \in Z^{n \times n}$ and $b \in \mathbb{R}^{n}$ by
\[
a_{ij}:=
\left\{
\begin{array}{cl}
\hat m_{ij} & \, \mbox{ if } i \not\in \Re, \\[1ex]
0 & \, \mbox{ if } i \in \Re \mbox{ and } i\not=j, \\[1ex]
\hat m_{ii}  & \, \mbox{ if } i \in \Re \mbox{ and } i=j,
\end{array}
\right.
\quad
\mbox{and}
\quad
b_{i}:=
\left\{
\begin{array}{cl}
\hat q_{i} & \, \mbox{ if } i \not\in \Re, \\[1ex]
0 & \, \mbox{ if } i \in \Re .
\end{array}
\right.
\]
For example, let $n=5$ and $\Re =\{2,4\}$. Then,
\[
A=
\left(
\begin{array}{ccccc}
\hat m_{11} & \hat m_{12} & \hat m_{13} & \hat m_{14} & \hat m_{15} \\[1ex]
0 & \hat m_{22} & 0 & 0 & 0 \\[1ex]
\hat m_{31} & \hat m_{32} & \hat m_{33} & \hat m_{34} & \hat m_{35} \\[1ex]
0 & 0 & 0 & \hat m_{44} & 0 \\[1ex]
\hat m_{51} & \hat m_{52} & \hat m_{53} & \hat m_{54} & \hat m_{55} 
\end{array}
\right)
\quad \mbox{and} \quad
b=
\left(
\begin{array}{c}
\hat q_{1} \\[1ex]
0  \\[1ex]
\hat q_{3}\\[1ex]
0 \\[1ex]
\hat q_{5}  
\end{array}
\right).
\]
We will show
\begin{equation} \label{zielkorfan}
b+A \tilde z \leq o \leq b+A \hat z.
\end{equation}
On the one hand, if $i \in \Re$ we have $\tilde z_i =0$ and
\[
\big( b + A \tilde z \big)_i =0 \leq \hat m _{ii} \hat z_i = \big( b + A \hat z \big)_i .
\]
On the other hand, if $i \not\in \Re$ we have $\tilde z_i >0$, $\big(\tilde q + \tilde M \tilde z \big)_i=0$ and therefore
\[
\begin{array}{rcl}
\big( b + A \tilde z\big)_i & = & \hat q_i + \sum \limits_{j=1}^{n} \hat m_{ij} \tilde z_j\\[2ex]
                                       & \leq & \tilde q_i + \sum \limits_{j=1}^{n} \tilde m_{ij} \tilde z_j=0\\[2ex]
                                       & \leq & \big(    \hat q + \hat M \hat z  \big)_i \\[2ex]
                                       & = & \big( b + A \hat z\big)_i .
\end{array}
\]
So we have shown (\ref{zielkorfan}) and  by Corollary \ref{korfan} $A$ is an M-matrix. Hence, we have $\tilde z \leq \hat z$.\\[1ex]
Case 2: Let $\tilde M$ be an M-matrix.\\[1ex]
Then $\tilde z $ is the unique solution of  the LCP defined by $\tilde M$ and $\tilde q$. Let $\Re =\{i: \tilde z_i=0\}$. Then we define $A \in Z^{n \times n}$ and $b \in \mathbb{R}^{n}$ by
\[
a_{ij}:=
\left\{
\begin{array}{cl}
\tilde m_{ij} & \, \mbox{ if } i \not\in \Re, \\[1ex]
0 & \, \mbox{ if } i \in \Re \mbox{ and } i\not=j, \\[1ex]
\tilde m_{ii}  & \, \mbox{ if } i \in \Re \mbox{ and } i=j,
\end{array}
\right.
\quad
\mbox{and}
\quad
b_{i}:=
\left\{
\begin{array}{cl}
\tilde q_{i} & \, \mbox{ if } i \not\in \Re, \\[1ex]
0 & \, \mbox{ if } i \in \Re .
\end{array}
\right.
\]
For example, let $n=5$ and $\Re =\{2,4\}$. Then,
\[
A=
\left(
\begin{array}{ccccc}
\tilde m_{11} & \tilde m_{12} & \tilde m_{13} & \tilde m_{14} & \tilde m_{15} \\[1ex]
0 & \tilde m_{22} & 0 & 0 & 0 \\[1ex]
\tilde m_{31} & \tilde m_{32} & \tilde m_{33} & \tilde m_{34} & \tilde m_{35} \\[1ex]
0 & 0 & 0 & \tilde m_{44} & 0 \\[1ex]
\tilde m_{51} & \tilde m_{52} & \tilde m_{53} & \tilde m_{54} & \tilde m_{55} 
\end{array}
\right)
\quad \mbox{and} \quad
b=
\left(
\begin{array}{c}
\tilde q_{1} \\[1ex]
0  \\[1ex]
\tilde q_{3}\\[1ex]
0 \\[1ex]
\tilde q_{5}  
\end{array}
\right).
\]
We will show
\begin{equation} \label{zielkorfan2}
b+A \tilde z \leq o \leq b+A \hat z.
\end{equation}
On the one hand, if $i \in \Re$ we have $\tilde z_i =0$ and
\[
\big( b + A \tilde z \big)_i =0 \leq \tilde m _{ii} \hat z_i = \big( b + A \hat z \big)_i .
\]
On the other hand, if $i \not\in \Re$ we have $\tilde z_i >0$, $\big(\tilde q + \tilde M \tilde z \big)_i=0$ and therefore
\[
\begin{array}{rcl}
\big( b + A \tilde z\big)_i & = & \tilde q_i + \sum \limits_{j=1}^{n} \tilde m_{ij} \tilde z_j=0\\[2ex]
                                       & \leq & \big(    \hat q + \hat M \hat z  \big)_i \\[2ex]
                                       & = & \hat q_i + \sum \limits_{j=1}^{n} \hat m_{ij} \hat z_j\\[2ex]
                                       & \leq & \tilde q_i + \sum \limits_{j=1}^{n} \tilde m_{ij} \hat z_j\\[2ex]
                                       & = & \big( b + A \hat z\big)_i .
\end{array}
\]
So we have shown (\ref{zielkorfan2}) and  by Corollary \ref{korfan} $A$ is an M-matrix. Hence, we have $\tilde z \leq \hat z$.\\[1ex]
Now we prove the supplement.\\[1ex]
Case 1: Let $\hat M^{-1} >O$\\[1ex]
Then we have
\[
\hat q + \hat M \tilde z \leq  \tilde q + \tilde M \tilde z =o \leq  \hat q + \hat M \hat z
\]
and
\[
\Big(
\hat q + \hat M \tilde z
\Big)_{i_0} < 
\Big(
\hat q + \hat M \hat z
\Big)_{i_0}.
\]
Hence, 
\[
o\leq \hat M (\hat z - \tilde z)
\]
and
\[
0< \Big(\hat M (  \hat z - \tilde z)
\Big)_{i_0}.
\]
With $\hat M^{-1} >O$ we get
\[
o< \hat M^{-1} \Big(
\hat M (\hat z - \tilde z)
\Big)= \hat z - \tilde z;
\]
i.e., $ \tilde z < \hat z$.\\[1ex]
Case 2: Let $\tilde  M^{-1} >O$\\[1ex]
Then we have
\[
\tilde q + \tilde M \tilde z =o \leq  \hat q + \hat M \hat z \leq   \tilde q + \tilde M \hat z
\]
and
\[
\Big(
\tilde q + \tilde M \tilde z
\Big)_{i_0} < 
\Big(
\tilde q + \tilde M \hat z
\Big)_{i_0}.
\]
Hence, 
\[
o\leq \tilde M (\hat z - \tilde z)
\]
and
\[
0< \Big(\tilde M (  \hat z - \tilde z)
\Big)_{i_0}.
\]
With $\tilde M^{-1} >O$ we get
\[
o< \tilde M^{-1} \Big(
\tilde M (\hat z - \tilde z)
\Big)= \hat z - \tilde z;
\]
i.e., $ \tilde z < \hat z$.\hfill $\Box$ 
\begin{example}
{\rm Let $\hat M = 0$, $\tilde M=1$, $\hat q=0$, and $\tilde q=3$. Then for any $t\geq 0$ $\hat z=t$ is a solution of the LCP defined by 
$\hat M $ and $\hat q$. On the other hand, $\tilde z =0$ is the unique solution of the LCP defined by 
$\tilde M $ and $\tilde q$. As predicted in Theorem \ref{neuestheorem} we have $\tilde z \leq \hat z$, since $\tilde M$ is an $1 \times 1$ M-matrix. The supplement of  Theorem \ref{neuestheorem}, however, cannot be applied, since neither $\tilde z >0$ nor $\tilde q= - \tilde M \tilde z$ is valid.}
\end{example}
\begin{example}
{\rm Let
\[
\hat M=
\left(
\begin{array}{cc}
2 & 0 \\
0 & 5
\end{array}
\right), \quad
\tilde M=
\left(
\begin{array}{cc}
2 & 7 \\
6 & 5
\end{array}
\right), \quad
\hat q=
\left(
\begin{array}{c}
-4 \\
-5
\end{array}
\right)=\tilde q.
\]
The unique solution of the LCP defined by  $\hat M $ and $\hat q$ is
\[
\hat z=
\left(
\begin{array}{c}
2 \\
1
\end{array}
\right).
\] 
The  LCP defined by  $\tilde M $ and $\tilde q$ has exactly three solutions
\[
\tilde z=
\frac{1}{32}
\left(
\begin{array}{c}
15 \\
14
\end{array}
\right),
\quad
\tilde z=
\left(
\begin{array}{c}
2 \\
0
\end{array}
\right),
\quad
\tilde z=
\left(
\begin{array}{c}
0 \\
1
\end{array}
\right).
\]
In all three cases we have $\tilde z \leq \hat z$ as predicted in Theorem \ref{neuestheorem} since $\hat M$ is an M-matrix. The supplement of  Theorem \ref{neuestheorem}, however, cannot be applied, since neither $\hat M^{-1} >O$ nor $\tilde M^{-1}>O$ is valid.}
\end{example}

\section{The solution set $\Sigma([M],[q])$}
We consider the $n \times 2n$  system of linear interval equations
\begin{equation} \label{intsys}
\Big( I \, \vdots \, -M \Big) \left(
\begin{array}{c}
w \\
z
\end{array}
\right)
 =q, \quad M \in [M], \, q \in [q]
\end{equation}
where $I$ denotes the identity. Due to Theorem 3.4.3 in \cite{neumaier} 
\[
\left(
\begin{array}{c}
w \\
z
\end{array}
\right) \in  \mathbb{R}^{2n}
\]
is a solution of (\ref{intsys}) if and only if
\[
\Big( I \, \vdots \, -[M] \Big) \left(
\begin{array}{c}
w \\
z
\end{array}
\right) \cap [q] \not= \emptyset;
\]
i.e., for $i=1,...,n$ we have
\[
\underline{w_i-\sum \limits_{j=1}^{n}[m_{ij}] \cdot z_j} \leq \overline{q}_i
\quad \mbox{and} \quad 
 \underline{q}_i \leq
\overline{w_i-\sum \limits_{j=1}^{n}[m_{ij}] \cdot z_j}. 
\]
Since $z \geq o$ these conditions simplify to
\[
w_i-\sum \limits_{j=1}^{n} \overline{m}_{ij} \cdot z_j \leq \overline{q}_i
\quad \mbox{and} \quad 
 \underline{q}_i \leq
w_i-\sum \limits_{j=1}^{n}\underline{m}_{ij} \cdot z_j, \quad i=1,...,n.
\]
In addition, it has to hold that $w^{\mbox{\scriptsize T}}z=0$, $w\geq o,\, z\geq o$; i.e.,
\[
w_i =0 \quad \mbox{or} \quad z_i=0 \quad \mbox{ for i=1,...,n.}
\]
So, we have to consider $2^n$ cases, and  taking into account $w\geq o$ and $z\geq o$ each case can give a contribution to the solution set or not.

\section{The symmetric solution set $\Sigma_{sym}([M],[q])$}
\subsection{The two-dimensional case}
The process described in Section 4 gives four cases:
\[
z_1=0, \, z_2=0, \quad 
\left(
\begin{array}{cc}
1 & 0 \\
0 & 1
\end{array}
\right)
 \left(
\begin{array}{c}
w_1 \\
w_2
\end{array}
\right) 
=q, \quad q \in [q],
\]
\[
w_1=0, \, z_2=0, \quad 
\left(
\begin{array}{cc}
-m_{11} & 0 \\
-m_{21} & 1
\end{array}
\right)
 \left(
\begin{array}{c}
z_1 \\
w_2
\end{array}
\right) 
=q, \quad q \in [q], \quad m_{11} \in [m_{11}], \, m_{21} \in [m_{21}],
\]
\[
z_1=0, \, w_2=0, \quad 
\left(
\begin{array}{cc}
1 & -m_{12} \\
0 & -m_{22}
\end{array}
\right)
 \left(
\begin{array}{c}
w_1 \\
z_2
\end{array}
\right) 
=q, \quad q \in [q], \quad m_{12} \in [m_{12}], \, m_{22} \in [m_{22}],
\]
\[
w_1=0, \, w_2=0, \quad 
\left(
\begin{array}{cc}
-m_{11} & -m_{12} \\
-m_{21} & -m_{22}
\end{array}
\right)
 \left(
\begin{array}{c}
z_1 \\
z_2
\end{array}
\right) 
=q, \quad q \in [q], \quad M \in [M].
\]
Using the symmetric structure of $M$ is only possible in the fourth case via $m_{21}=m_{12}$. Then, we have
\[
\begin{array}{ccccc}
\underline{q}_{1} & \leq & -m_{11} \cdot z_1 -m_{12} \cdot z_2 & \leq & \overline{q}_1 \\[1ex]
\underline{q}_{2} & \leq & -m_{12} \cdot z_1 -m_{22} \cdot z_2 & \leq & \overline{q}_2 
\end{array}
\]
which leads to
\[
\begin{array}{ccccc}
\underline{q}_{1} + \underline{m}_{11} \cdot z_1 & \leq &  -m_{12} \cdot z_2 & \leq & \overline{q}_1 + \overline{m}_{11} z_1 \\[1ex]
\underline{q}_{2} + \underline{m}_{22} \cdot z_2 & \leq & -m_{12} \cdot z_1  & \leq & \overline{q}_2 + \overline{m}_{22} z_2
\end{array}
\]
which leads to
\[
\begin{array}{ccccc}
\displaystyle{\frac{\underline{q}_{1} + \underline{m}_{11} \cdot z_1}{z_2}} & \leq &  -m_{12} & \leq & \displaystyle{\frac{\overline{q}_1 + \overline{m}_{11} z_1}{z_2}} \\[3ex]
\displaystyle{\frac{\underline{q}_{2} + \underline{m}_{22} \cdot z_2}{z_1}} & \leq &  -m_{12} & \leq & \displaystyle{\frac{\overline{q}_2 + \overline{m}_{22} z_2}{z_1}}
\end{array}
\]
which leads to
\[
\begin{array}{ccc}
(\underline{q}_{1} + \underline{m}_{11} \cdot z_1)\cdot z_1 & \leq & (\overline{q}_2 + \overline{m}_{22} z_2)\cdot z_2 \\[1ex]
(\underline{q}_{2} + \underline{m}_{22} \cdot z_2)\cdot z_2 & \leq & (\overline{q}_1 + \overline{m}_{11} z_1)\cdot z_1
\end{array}
\]
which results in
\[
\begin{array}{ccc}
0 & \leq &  \overline{m}_{22} z_2^2   - \underline{m}_{11} \cdot z_1^2 + \overline{q}_2 \cdot z_2 -  \underline{q}_{1} \cdot z_1    \\[1ex]
0 & \leq &  \overline{m}_{11} z_1^2   - \underline{m}_{22} \cdot z_2^2 + \overline{q}_1 \cdot z_1 -  \underline{q}_{2} \cdot z_2.  
\end{array}
\]

\subsection{The three-dimensional case}
\subsubsection{The Fourier-Motzkin elimination process for the case that $w_1=0$, $w_2=0$, and $w_3=0$}
Case 1: $m_{12}=m_{21}$. Then,
\[
\left\{
\begin{array}{ccccc}
\underline{q}_{1} & \leq & -m_{11} \cdot z_1-m_{12} \cdot z_2 - m_{13}\cdot z_3 & \leq & \overline{q}_{1}  \\[1ex]
\underline{q}_{2} & \leq & -m_{21} \cdot z_1-m_{22} \cdot z_2 - m_{23}\cdot z_3 & \leq & \overline{q}_{2}
\end{array}
\right.
\] 
which is
\[
\left\{
\begin{array}{ccccc}
\underline{q}_{1} + m_{11} \cdot z_1 +  m_{13}\cdot z_3 & \leq & -m_{12} \cdot z_2  & \leq & \overline{q}_{1} + m_{11} \cdot z_1 +  m_{13}\cdot z_3  \\[1ex]
\underline{q}_{2} + m_{22} \cdot z_2 + m_{23}\cdot z_3 & \leq & -m_{21} \cdot z_1 & \leq & \overline{q}_{2} + m_{22} \cdot z_2 + m_{23}\cdot z_3
\end{array}
\right.
\] 
and
\[
\left\{
\begin{array}{ccccc}
\displaystyle{\frac{\underline{q}_{1} + \underline{m}_{11} \cdot z_1 +  \underline{m}_{13}\cdot z_3}{z_2} }& \leq & -m_{12}   & \leq & \displaystyle{\frac{\overline{q}_{1} +\overline{ m}_{11} \cdot z_1 + \overline{ m}_{13}\cdot z_3}{z_2} } \\[3ex]
\displaystyle{\frac{\underline{q}_{2} + \underline{m}_{22} \cdot z_2 + \underline{m}_{23}\cdot z_3}{z_1}} & \leq & -m_{21}  & \leq & \displaystyle{\frac{\overline{q}_{2} + \overline{m}_{22} \cdot z_2 + \overline{m}_{23}\cdot z_3}{z_1}}
\end{array}
\right.
\] 
and, since  $m_{12}=m_{21}$,
\[
\left\{
\begin{array}{ccc}
(\underline{q}_{1} +\underline{ m}_{11} \cdot z_1 +  \underline{m}_{13}\cdot z_3 ) \cdot z_1&  \leq & (\overline{q}_{2} + \overline{m}_{22} \cdot z_2 + \overline{m}_{23}\cdot z_3)\cdot z_2  \\[1ex]
(\underline{q}_{2} + \underline{m}_{22} \cdot z_2 + \underline{m}_{23}\cdot z_3) \cdot z_2 & \leq & (\overline{q}_{1} +\overline{ m}_{11} \cdot z_1 + \overline{ m}_{13}\cdot z_3 ) \cdot z_1
\end{array}
\right.
\] 
which finally leads to
\[
\left\{
\begin{array}{ccc}
0 & \leq & \overline{m}_{22} \cdot z_2^{2} + \overline{m}_{23}\cdot z_2 z_3 +\overline{q}_{2} \cdot z_2 - 
\underline{m}_{11} \cdot z_1^2 - \underline{m}_{13} \cdot z_1 z_3  -\underline{q}_{1}  \cdot z_1  \\[1ex]
0 & \leq & \overline{m}_{11} \cdot z_1^{2} + \overline{m}_{13}\cdot z_1 z_3 +\overline{q}_{1} \cdot z_1 - 
\underline{m}_{22} \cdot z_2^2 - \underline{m}_{23} \cdot z_2 z_3  -\underline{q}_{2}  \cdot z_2
\end{array}
\right.
\] 
Case 2: $m_{13}=m_{31}$. Analogously to Case 1 we get
\[
\left\{
\begin{array}{ccc}
0 & \leq & \overline{m}_{33} \cdot z_3^{2} + \overline{m}_{32}\cdot z_2 z_3 +\overline{q}_{3} \cdot z_3 - 
\underline{m}_{11} \cdot z_1^2 - \underline{m}_{12} \cdot z_1 z_2  -\underline{q}_{1}  \cdot z_1  \\[1ex]
0 & \leq & \overline{m}_{11} \cdot z_1^{2} + \overline{m}_{12}\cdot z_1 z_2 +\overline{q}_{1} \cdot z_1 - 
\underline{m}_{33} \cdot z_3^2 - \underline{m}_{32} \cdot z_2 z_3  -\underline{q}_{3}  \cdot z_3
\end{array}
\right.
\] 
Case 3: $m_{23}=m_{32}$. Analogously to Case 1 we get
\[
\left\{
\begin{array}{ccc}
0 & \leq & \overline{m}_{33} \cdot z_3^{2} + \overline{m}_{31}\cdot z_1 z_3 +\overline{q}_{3} \cdot z_3 - 
\underline{m}_{22} \cdot z_2^2 - \underline{m}_{21} \cdot z_1 z_2  -\underline{q}_{2}  \cdot z_2  \\[1ex]
0 & \leq & \overline{m}_{22} \cdot z_2^{2} + \overline{m}_{21}\cdot z_1 z_2 +\overline{q}_{2} \cdot z_2 - 
\underline{m}_{33} \cdot z_3^2 - \underline{m}_{31} \cdot z_1 z_3  -\underline{q}_{3}  \cdot z_3
\end{array}
\right.
\] 
\begin{theorem}
Let any $M \in [M]$ be a  P-matrix, then each quadric that appears during the Fourier-Motzkin elimination process for the case $w_1=0$, $w_2=0$, and $w_3=0$ is neither an ellipsoid neither an elliptic paraboloid neither an elliptic cylinder nor a parabolic cylinder. 
\end{theorem}
{\it Proof:} Case 1: $m_{12}=m_{21}$. Then,
\[
\left\{
\begin{array}{rcl}
0 & = & \overline{m}_{22} \cdot z_2^{2} + \overline{m}_{23}\cdot z_2 z_3 +\overline{q}_{2} \cdot z_2 - 
\underline{m}_{11} \cdot z_1^2 - \underline{m}_{13} \cdot z_1 z_3  -\underline{q}_{1}  \cdot z_1  \\[3ex]
& = & ( z_1 \, \, z_2 \, \, z_3) 
\left(
\begin{array}{ccc}
- \underline{m}_{11} & 0 & - \frac{1}{2} \underline{m}_{13} \\[1ex]
0 & \overline{m}_{22} & \frac{1}{2} \overline{m}_{23}\\[1ex]
 - \frac{1}{2} \underline{m}_{13} &  \frac{1}{2} \overline{m}_{23} & 0
\end{array}
\right) 
\cdot 
\left(
\begin{array}{c}
z_1 \\
z_2\\
z_3
\end{array}
\right)
+ 
 ( -\underline{q}_1 \, \, \overline{q}_2 \, \, 0)
 \cdot 
\left(
\begin{array}{c}
z_1 \\
z_2\\
z_3
\end{array}
\right) .
\end{array}
\right.
\] 
We consider
\[
p(\lambda) =
\det \left(
\begin{array}{ccc}
- \underline{m}_{11} - \lambda & 0 & - \frac{1}{2} \underline{m}_{13} \\[1ex]
0 & \overline{m}_{22} - \lambda & \frac{1}{2} \overline{m}_{23}\\[1ex]
 - \frac{1}{2} \underline{m}_{13} &  \frac{1}{2} \overline{m}_{23} & -\lambda
\end{array}
\right) 
\]
\[
= ( \underline{m}_{11} + \lambda )(\overline{m}_{22} - \lambda) \lambda
-  \frac{1}{4} \underline{m}_{13} ^2 ( \overline{m}_{22} - \lambda) + ( \underline{m}_{11} + \lambda)\frac{1}{4} \overline{m}_{23} ^2 
\]
Subcase 1: $\underline{m}_{13} \cdot \overline{m}_{23} \not = 0$. Then,
\[
p(-\underline{m}_{11} )< 0, \quad p(\overline{m}_{22}) > 0,  \quad \lim \limits_{\lambda \to -\infty}p(\lambda)=\infty, \quad \mbox{and} \quad
\lim \limits_{\lambda \to \infty}p(\lambda)=-\infty.
\]
So, there are three different eigenvalues and they don't have the same sign. Even if one eigenvalue is 0, then the other two eigenvalues have different signs. Therefore, the quadric cannot  be an ellipsoid, it cannot be  an elliptic paraboloid, it cannot be  an elliptic cylinder, and it cannot be a parabolic cylinder.

Subcase 2: $\underline{m}_{13}=0, \,  \overline{m}_{23} \not = 0$. Then,
\[
p(\lambda)=
 -( \underline{m}_{11} + \lambda )(\lambda^2-\overline{m}_{22}  \lambda
-  \frac{1}{4} \overline{m}_{23} ^2) 
\]
and
the eigenvalues are
\[
\lambda_1=-\underline{m}_{11}<0
\]
and 
\[
\lambda_{2,3}=\frac{\overline{m}_{22}}{2}\pm
\sqrt{
\left(
\frac{\overline{m}_{22}}{2}
\right)^2
+
\frac{1}{4} \overline{m}_{23} ^2}
\]
with
\[
\lambda_{2}<0, \quad \lambda_{3}>0.
\]
So, the eigenvalues  don't have the same sign. Therefore, the quadric cannot be  an ellipsoid, it cannot be  an elliptic paraboloid, it cannot be a hyperbolic paraboloid, it cannot be  an elliptic cylinder, it cannot be a hyperbolic cylinder, and it cannot be a parabolic cylinder.

Subcase 3: $\underline{m}_{13}\not=0, \,  \overline{m}_{23}  = 0$. This subcase runs analogously to subcase 2 and is omitted here.

Subcase 4: $\underline{m}_{13}=0, \,  \overline{m}_{23}  = 0$. Then
\[
0= \overline{m}_{22} \cdot z_2^{2}  +\overline{q}_{2} \cdot z_2 - 
\underline{m}_{11} \cdot z_1^2   -\underline{q}_{1}  \cdot z_1.
\]
This is a hyperbolic cylinder or (as a special case of it) two intersecting planes.

Case 2: $m_{13}=m_{31}$ and Case 3: $m_{32}=m_{23}$ run analogously and are omitted here.

\subsubsection{The Fourier-Motzkin elimination process for the case that exactly two $w_i=0$ and exactly one $z_j=0$}
Without loss of generality let $w_1=0$, $w_2=0$, and $z_3=0$. Then we have to solve the system of linear equations
\[
\left(
\begin{array}{ccc}
-m_{11} & - m_{12} & 0 \\[1ex]
-m_{21} & - m_{22} & 0 \\[1ex]
 -m_{31} & - m_{32} & 1
\end{array}
\right)
\left(
\begin{array}{c}
z_1\\[1ex]
z_2 \\[1ex]
w_3
\end{array}
\right)=
\left(
\begin{array}{c}
q_1\\[1ex]
q_2 \\[1ex]
q_3
\end{array}
\right).
\]
The Fourier-Motzkin elimination process can only be applied for the case $m_{12}=m_{21}$. I.e.,
\[
\left\{
\begin{array}{ccccc}
\underline{q}_{1} & \leq & -m_{11} \cdot z_1-m_{12} \cdot z_2   & \leq & \overline{q}_{1}  \\[1ex]
\underline{q}_{2} & \leq & -m_{21} \cdot z_1-m_{22} \cdot z_2  & \leq & \overline{q}_{2}
\end{array}
\right.
\] 
which is
\[
\left\{
\begin{array}{ccccc}
\underline{q}_{1} + m_{11} \cdot z_1  & \leq & -m_{12} \cdot z_2  & \leq & \overline{q}_{1} + m_{11} \cdot z_1   \\[1ex]
\underline{q}_{2} + m_{22} \cdot z_2 & \leq & -m_{21} \cdot z_1 & \leq & \overline{q}_{2} + m_{22} \cdot z_2 
\end{array}
\right.
\] 
and
\[
\left\{
\begin{array}{ccccc}
\displaystyle{\frac{\underline{q}_{1} + \underline{m}_{11} \cdot z_1 }{z_2} }& \leq & -m_{12}   & \leq & \displaystyle{\frac{\overline{q}_{1} +\overline{ m}_{11} \cdot z_1 }{z_2} } \\[3ex]
\displaystyle{\frac{\underline{q}_{2} + \underline{m}_{22} \cdot z_2 }{z_1}} & \leq & -m_{21}  & \leq & \displaystyle{\frac{\overline{q}_{2} + \overline{m}_{22} \cdot z_2 }{z_1}}
\end{array}
\right.
\] 
and, since  $m_{12}=m_{21}$,
\[
\left\{
\begin{array}{ccc}
(\underline{q}_{1} +\underline{ m}_{11} \cdot z_1  ) \cdot z_1&  \leq & (\overline{q}_{2} + \overline{m}_{22} \cdot z_2 )\cdot z_2  \\[1ex]
(\underline{q}_{2} + \underline{m}_{22} \cdot z_2 ) \cdot z_2 & \leq & (\overline{q}_{1} +\overline{ m}_{11} \cdot z_1  ) \cdot z_1
\end{array}
\right.
\] 
which finally leads to
\[
\left\{
\begin{array}{ccc}
0 & \leq & \overline{m}_{22} \cdot z_2^{2}  +\overline{q}_{2} \cdot z_2 - 
\underline{m}_{11} \cdot z_1^2  -\underline{q}_{1}  \cdot z_1  \\[1ex]
0 & \leq & \overline{m}_{11} \cdot z_1^{2}  +\overline{q}_{1} \cdot z_1 - 
\underline{m}_{22} \cdot z_2^2   -\underline{q}_{2}  \cdot z_2 .
\end{array}
\right.
\] 
So, immediatly we get the following theorem.
\begin{theorem}
Let any $M \in [M]$ be a  P-matrix, then each quadric that appears during the Fourier-Motzkin elimination process for the case that exactly two  $w_i=0$ and exactly one  $z_j=0$ is 
a hyperbolic cylinder or (as a special case of it) two intersecting planes.
\end{theorem}

\section{Examples}

\subsection{Two-dimensional examples}
\subsubsection{All $M \in [M]$ are M-matrices}
Let
\[
[M]=\left(
\begin{array}{cc}
[\frac{1}{8},1] & [-\frac{1}{4},-\frac{1}{5}] \\[1ex]
[-\frac{1}{4},-\frac{1}{5}] & 1
\end{array}
\right) \quad \mbox{and} \quad
[q]=\left(
\begin{array}{c}
[-3,-1]  \\[1ex]
[1,2] 
\end{array}
\right).
\] 
All $M \in [M]$ are M-matrices due to 
$\underline{M} 
\left(
\begin{array}{c}
3 \\ 1
\end{array}
\right)
> o$, Lemma \ref{lemmafan}, and Corollary \ref{korfan}. The process described in Section 4 gives 
\[
\begin{array}{rclcrcl}
w_1 -z_1+\frac{1}{5}z_2 & \leq & -1 & \quad \mbox{and} \quad & -3 & \leq & w_1 - \frac{1}{8}z_1 + \frac{1}{4} z_2 \\[1ex]
w_2 +\frac{1}{5}z_1 -z_2 & \leq & 2 & \quad \mbox{and} \quad & 1 & \leq & w_2 + \frac{1}{4} z_1-z_2
\end{array}
\]

Case 1: $w_1 \geq 0, \, w_2\geq 0,\, z_1=0,\, z_2=0$. We get 
the contradiction $w_1\in [-3,-1]$
 and $w_1\geq 0$. Therefore, Case 1 gives no contribution to the solution set.

Case 2: $w_1 \geq 0, \, z_2\geq 0,\, z_1=0,\, w_2=0$. We get 
the contradiction $1\leq -z_2$
 and $z_2\geq 0$. Therefore, Case 2 gives no contribution to the solution set.

Case 3: $z_1 \geq 0, \, w_2\geq 0,\, w_1=0,\, z_2=0$. The four inequalities lead to
\[
\begin{array}{rclcrcl}
1 & \leq & z_1 & \quad \mbox{and} \quad & z_1 & \leq & 24 \\[1ex]
w_2 & \leq & 2- \frac{1}{5}z_1  & \quad \mbox{and} \quad & 1-\frac{1}{4} z_1 & \leq & w_2 
\end{array}
\]

Via $z_1 \geq 0, \, w_2\geq 0$ we get the polygon depicted in the following figure.
\begin{figure}[h]
\unitlength1.0cm
\begin{picture}(12,4)
\thinlines
\put(1,0){\vector(0,1){3.8}}
\put(0,1){\vector(1,0){12.5}}
\put(2,0.9){\line(0,1){0.1}}
\put(3,0.9){\line(0,1){0.1}}
\put(4,0.9){\line(0,1){0.1}}
\put(5,0.9){\line(0,1){0.1}}
\put(6,0.9){\line(0,1){0.1}}
\put(7,0.9){\line(0,1){0.1}}
\put(8,0.9){\line(0,1){0.1}}
\put(9,0.9){\line(0,1){0.1}}
\put(10,0.9){\line(0,1){0.1}}
\put(11,0.9){\line(0,1){0.1}}
\put(0.9, 2){\line(1,0){0.1}}
\put(0.9, 3){\line(1,0){0.1}}
\put(0.7,1.9){$1$}
\put(0.7,2.9){$2$}
\put(0.35,3.5){$w_2$}
\put(5.9,0.5){$5$}
\put(10.8,0.5){$10$}
\put(12,0.5){$z_{1}$}
\thicklines
\put(2,2.8){\line(5,-1){9}}
\put(2,1.75){\line(4,-1){3}}
\put(5,1){\line(1,0){6}}
\put(2,2.8){\line(0,-1){1.05}}
\end{picture}
\end{figure}

The contribution to the solution set is
\[
\left\{
\left(
\begin{array}{c}
z_1 \\
z_2
\end{array}
\right) \in  \mathbb{R}^2\, : \, z_1 \in [1,10], \, z_2=0
\right\} .
\]

Case 4: $z_1 \geq 0, \, z_2\geq 0,\, w_1=0,\, w_2=0$. The four inequalities lead to
\[
\begin{array}{rclcrcl}
z_2 & \leq & 5z_1 -5 & \quad \mbox{and} \quad & \frac{1}{2}z_1-12 & \leq & z_2 \\[1ex]
\frac{1}{5}z_1 -2 & \leq & z_2  & \quad \mbox{and} \quad & z_2 & \leq & \frac{1}{4}z_1-1 
\end{array}
\]

Via $z_1 \geq 0, \, z_2\geq 0$ we get the polygon with its vertices 
\[
A(4|0), \,  B(10|0), \,  C(\frac{100}{3}|\frac{14}{3}), \mbox{ and }\, D(44|10).
\] 
It is depicted in the following figure.
\begin{figure}[h]
\unitlength1.0cm
\begin{picture}(12,5)
\thinlines
\put(1,0){\vector(0,1){4.8}}
\put(0,1){\vector(1,0){12.5}}
\put(3.5,0.9){\line(0,1){0.1}}
\put(2,0.9){\line(0,1){0.1}}
\put(1.9,0.5){$4$}
\put(6,0.9){\line(0,1){0.1}}
\put(8.5,0.9){\line(0,1){0.1}}
\put(11,0.9){\line(0,1){0.1}}
\put(0.9, 3.5){\line(1,0){0.1}}
\put(0.5,3.4){$10$}
\put(0.5,4.3){$z_2$}
\put(3.3,0.5){$10$}
\put(5.8,0.5){$20$}
\put(8.3,0.5){$30$}
\put(10.8,0.5){$40$}
\put(12,0.5){$z_{1}$}
\thicklines
\put(3.5,1){\line(5,1){5.83}}
\put(2,1){\line(4,1){10}}
\put(9.33,2.16){\line(2,1){2.66}}
\put(9,2.3){$L$}
\end{picture}
\end{figure}

So, the solution set is
\[
\Sigma([M],[q ])=\left(
\begin{array}{c}
[1,10] \\[1ex]
0
\end{array}
\right)
\cup
\left\{
\left(
\begin{array}{c}
z_1 \\[1ex]
z_2
\end{array}
\right) \in  \mathbb{R}^2\, : \, (z_1,z_2) \in L
\right\}
\]
where $L$ is depicted in the following figure. $\Sigma([M],[q])$ is connected, but is not convex with
\[
\inf \Sigma([M],[q]) =\left(
\begin{array}{c}
1 \\[1ex]
0
\end{array}
\right)
\quad
\mbox{and} \quad
\sup \Sigma([M],[q]) =\left(
\begin{array}{c}
44 \\[1ex]
10
\end{array}
\right).
\]

According to Theorem \ref{neuestheorem}, $\inf \Sigma([M],[q])$
is the unique solution of the LCP defined by
\[
\overline{M}=\left(
\begin{array}{cc}
1 & -\frac{1}{5} \\[1ex]
-\frac{1}{5} & 1
\end{array}
\right) \quad \mbox{and} \quad
\overline{q}=\left(
\begin{array}{c}
-1  \\[1ex]
2 
\end{array}
\right)
\] 
and $\sup \Sigma([M],[q])$
is the unique solution of the LCP defined by
\[
\underline{M}=\left(
\begin{array}{cc}
\frac{1}{8} & -\frac{1}{4} \\[1ex]
-\frac{1}{4} & 1
\end{array}
\right) \quad \mbox{and} \quad
\underline{q}=\left(
\begin{array}{c}
-3  \\[1ex]
1
\end{array}
\right).
\] 

Now, we consider the shape of the symmetric solution set. The process described in Section 5.1 gives  
\[
0\leq z_1^2 -z_1-z_2^2-z_2 \quad \mbox{and} \quad 
0 \leq z_2^2 +2 z_2-\frac{1}{8}z_1^2+3z_1
\]
which leads to
\[
0\leq (z_1-\frac{1}{2})^2 -(z_2+\frac{1}{2})^2 \quad \mbox{and} \quad 
0 \leq (z_2+1)^2 -\frac{1}{8}(z_1-12)^2+17.
\]

Using GeoGebra, for example, we see that the solution set is intersected by the hyperbola.\newline
\includegraphics[scale=0.24]{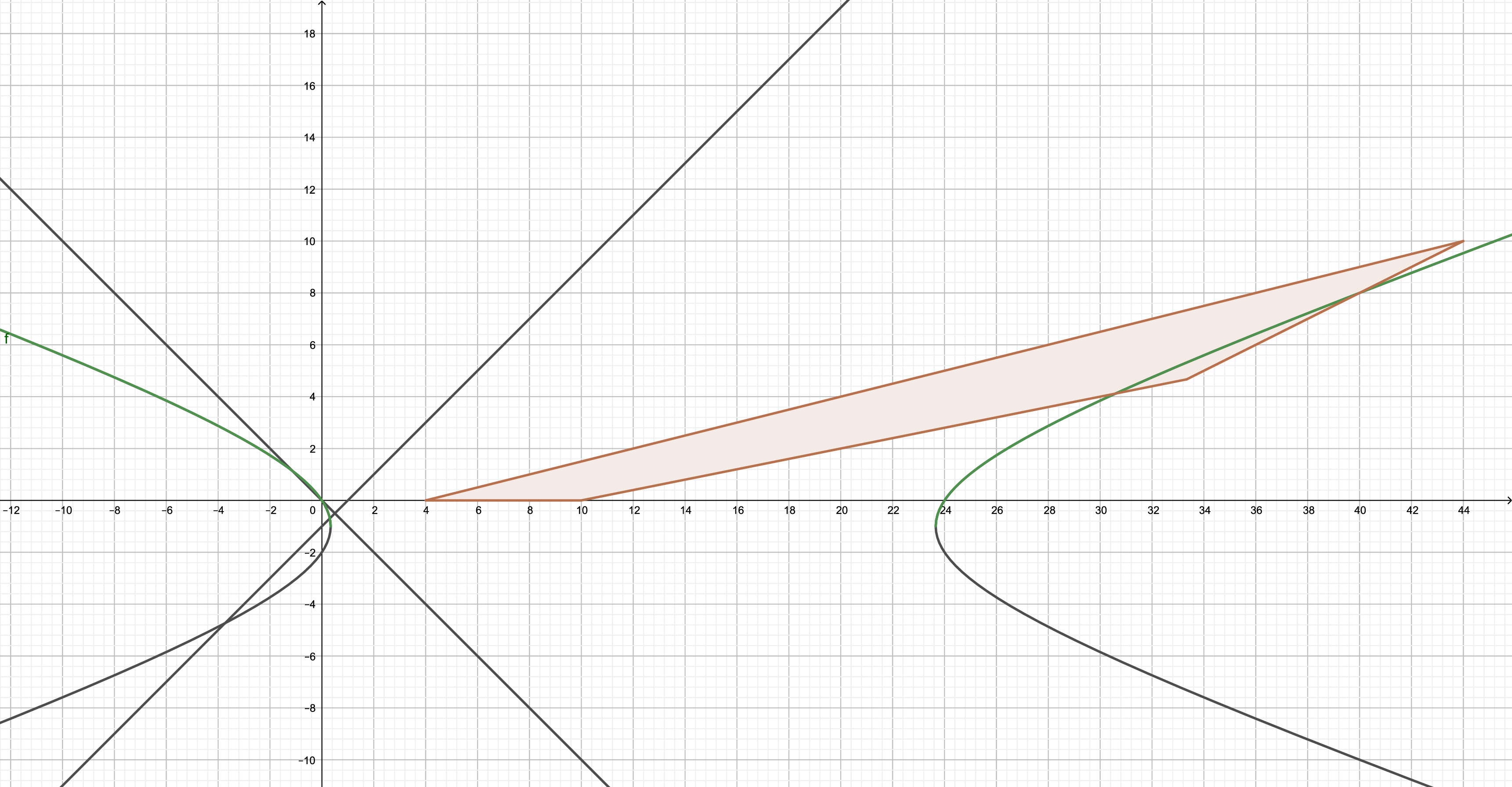}

So, the symmetric solution set is smaller than the solution set. The hyperbola cuts off a small region that contains the vertex $ C(\frac{100}{3}|\frac{14}{3})$; i.e., 
\[
z=
\left(
\begin{array}{c}
\displaystyle{\frac{100}{3}} \\[1ex]
\displaystyle{\frac{14}{3}}
\end{array}
\right)
\]
does not belong to the symmetric solution set, but it belongs to the solution set. It is the unique solution of the LCP defined by
\[
M=\left(
\begin{array}{cc}
\frac{1}{8} & -\frac{1}{4} \\[1ex]
-\frac{1}{5} & 1
\end{array}
\right) \quad \mbox{and} \quad
q=\left(
\begin{array}{c}
-3  \\[1ex]
2
\end{array}
\right).
\] 
Assume there is a symmetric matrix $M \in [M]$ and a $q \in [q]$ such that $z$ solves the LCP defined by $M$ and $q$. Then, since $z>o$, we have $q+Mz=o$. Hence, 
\[
q_1+m_{11} \cdot \frac{100}{3} + m_{12} \cdot \frac{14}{3}=0\quad \mbox{and} \quad
q_2+ m_{21}  \cdot \frac{100}{3}+  \frac{14}{3}=0.
\]
Since $m_{12}=m_{21}$ we get
\[
m_{11}=\frac{3 \cdot q_2 +14}{10000}\cdot 14 - q_1 \cdot \frac{3}{100} \in 
\frac{3 \cdot [1,2] +14}{10000}\cdot 14 - [-3,-1] \cdot \frac{3}{100}=[0.0538,0.118]
\]
which is a contradiction to $m_{11} \in [0.125,1]$.

\subsubsection{All $M \in [M]$ are $H_+$-matrices}
Let
\[
[M]=\left(
\begin{array}{cc}
[4,5] & [-1,2] \\[1ex]
[-1,2] & [2,3]
\end{array}
\right) \quad \mbox{and} \quad
[q]=\left(
\begin{array}{c}
[-2,-1]  \\[1ex]
[-1,1] 
\end{array}
\right).
\] 
By Lemma \ref{lemmafan} and Corollary \ref{korfan} every $M \in [M]$ is an $H_+$-matrix, since
\[
\left(
\begin{array}{rr}
4 & -2 \\[1ex]
-2 & 2
\end{array}
\right) \cdot \left(
\begin{array}{c}
1\\[1ex]
\displaystyle{\frac{3}{2}} 
\end{array}
\right)>o.
\]
The process described in Section 4 gives 
\[
\begin{array}{rclcrcl}
w_1 -5z_1 -2z_2 & \leq & -1 & \quad \mbox{and} \quad & -2 & \leq & w_1 - 4z_1 +  z_2 \\[1ex]
w_2 -2z_1 -3z_2 & \leq & 1 & \quad \mbox{and} \quad & -1 & \leq & w_2 + z_1-2z_2
\end{array}
\]

Case 1: $w_1 \geq 0, \, w_2\geq 0,\, z_1=0,\, z_2=0$. We get 
the contradiction $w_1 \leq -1$. Therefore, Case 1 gives no contribution to the solution set.

Case 2: $w_1 = 0, \, z_2= 0,\, z_1\geq 0,\, w_2\geq 0$. We get
\[
\begin{array}{rclcrcl}
\frac{1}{5} & \leq & z_1 & \quad \mbox{and} \quad & z_1 & \leq & \frac{1}{2} \\[1ex]
w_2 & \leq & 2z_1 + 1 & \quad \mbox{and} \quad & -z_1-1 & \leq & w_2 
\end{array}
\]
The contribution to the solution set is
\[
\left\{
\left(
\begin{array}{c}
z_1 \\
z_2
\end{array}
\right) \in  \mathbb{R}^2\, : \, z_1 \in [\frac{1}{5},\frac{1}{2}], \, z_2=0
\right\}.
\]
Case 3: $w_1 \geq 0, \, z_2\geq 0,\, z_1=0,\, w_2=0$. The four inequalities lead to
\[
\begin{array}{rclcrcl}
w_1 & \leq &2z_1 -1 & \quad \mbox{and} \quad & -z_2-2 & \leq & w_1 \\[1ex]
0 & \leq & z_2  & \quad \mbox{and} \quad & z_2 & \leq & \frac{1}{2} 
\end{array}
\]
The contribution to the solution set is
\[
\left\{
\left(
\begin{array}{c}
0 \\[1ex]
\frac{1}{2}
\end{array}
\right) 
\right\}.
\]
Case 4: $z_1 \geq 0, \, z_2\geq 0,\, w_1=0,\, w_2=0$. The four inequalities lead to
\[
\begin{array}{rclcrcl}
-\frac{5}{2}z_1 +\frac{1}{2}& \leq & z_2 & \quad \mbox{and} \quad & 4z_1-2 & \leq & z_2 \\[1ex]
-\frac{2}{3}z_1 -\frac{1}{3} & \leq & z_2  & \quad \mbox{and} \quad & z_2 & \leq & \frac{1}{2}z_1+\frac{1}{2} 
\end{array}
\]
So, the solution set is the polygon in the  figure below. Now, we consider
the shape of the symmetric solution set. The process described in Section 5.1 gives 
\[
0\leq 3z_2^2-4z_1^2 +z_2+2z_1= 3(z_2+\frac{1}{6})^2-4(z_1-\frac{1}{4})^2+\frac{1}{6} 
\]
and
\[
0\leq 5z_1^2-2z_2^2-z_1+z_2= 5(z_1-\frac{1}{10})^2 - 2(z_2-\frac{1}{4})^2+\frac{3}{40}.
\]
Using GeoGebra, for example, we see that the solution set is intersected by the hyperbola.\newline
\includegraphics[scale=0.5]{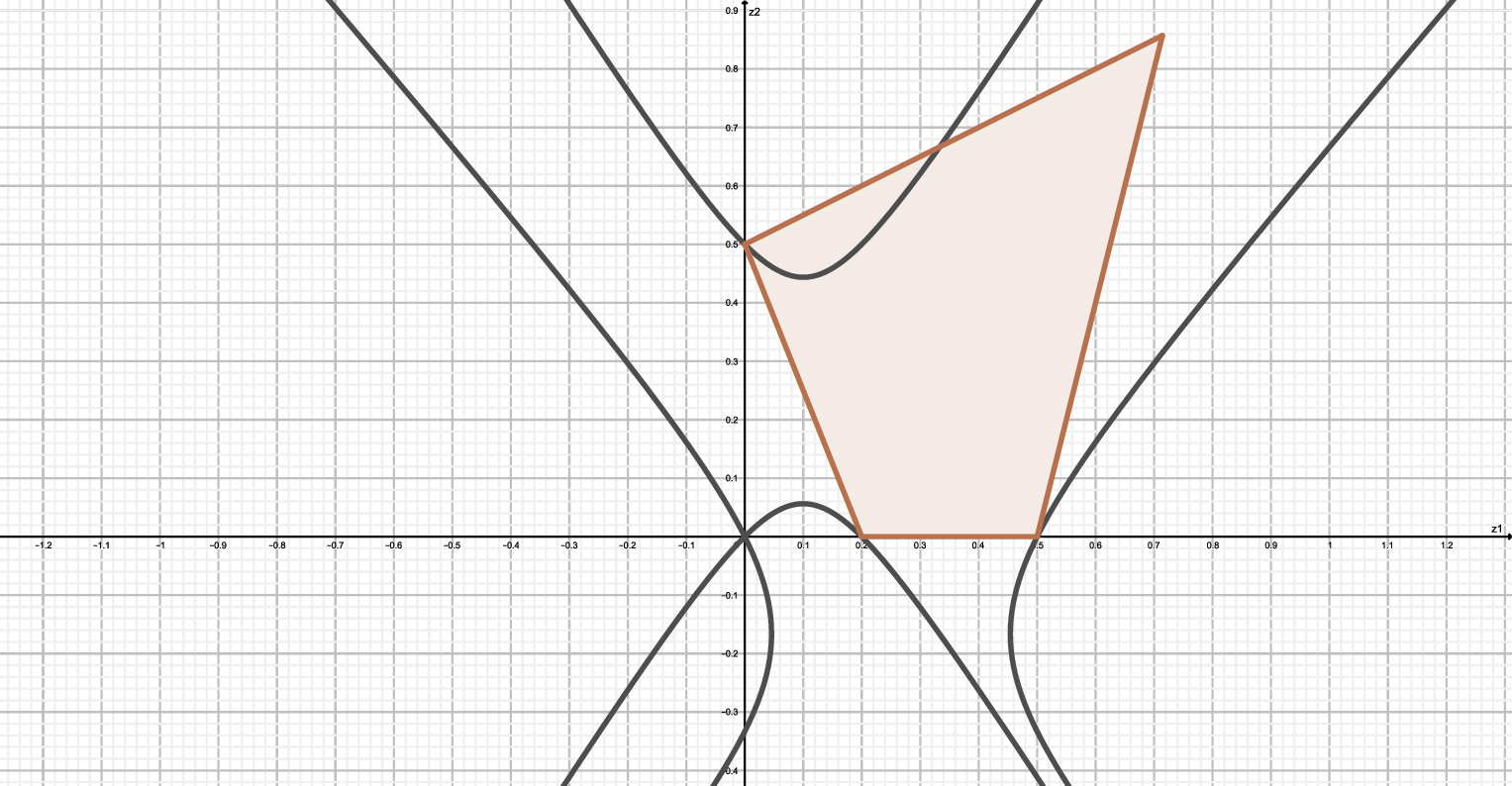}

Since $\underline{M}$ is an M-matrix we have due to Theorem \ref{neuestheorem} that $\underline{z}=\sup \Sigma([M],[q])$ where $\underline{z}$ is the unique solution of the LCP defined be $\underline{M}$ and $\underline{q}$. Furthermore, $\underline{M}$ is symmetric. So, $\underline{z}=\sup \Sigma_{sym}([M],[q])$. In addition, we have
\[
\underline{M}^{-1}=\frac{1}{7}
\left(
\begin{array}{rr}
2 & 1 \\[1ex]
1 & 4
\end{array}
\right) >O.
\]
By the second part of  Theorem \ref{neuestheorem} we have 
\[
z < \sup \Sigma_{sym}([M],[q]) \quad \mbox{for all } z \in \Sigma_{sym}([M],[q]) \backslash \{\underline{z} \}.
\]

\subsection{Three-dimensional examples}
\subsubsection{All $M \in [M]$ are P-matrices}
Let
\[
[M]=\left(
\begin{array}{ccc}
1 & [0,\frac{1}{2}] & [0,\frac{1}{2}] \\[1ex]
[0,\frac{1}{2}] & 1 & [0,\frac{1}{2}] \\[1ex]
[0,\frac{1}{2}] & [0,\frac{1}{2}] & 1
\end{array}
\right) \quad \mbox{and} \quad
[q]=\left(
\begin{array}{c}
-6\\[1ex]
[1,2] \\[1ex]
[-3,-2]
\end{array}
\right).
\] 
All $M \in [M]$ are P-matrices, see Example 3.1 in \cite{schaeferLAA}. The process described in Section 4 gives
\[
\begin{array}{rclcrcl}
w_1 -z_1-\frac{1}{2}z_2 -\frac{1}{2}z_3 & \leq & -6 & \quad \mbox{and} \quad & -6 & \leq & w_1 - z_1 \\[1ex]
w_2 -\frac{1}{2}z_1 -z_2 -\frac{1}{2}z_3 & \leq & 2 & \quad \mbox{and} \quad & 1 & \leq & w_2 -z_2\\[1ex]
w_3 -\frac{1}{2}z_1 -\frac{1}{2}z_2 -z_3 & \leq & -2 & \quad \mbox{and} \quad & -3 & \leq & w_3 -z_3.
\end{array}
\]

Case 1: $w_1\geq 0$, $w_2 \geq 0$, $w_3 \geq 0$, $z_1=0$, $z_2=0$, $z_3=0$. This leads to $-6=w_1\geq 0$, which is a contradiction. So, we have  no contribution to the solution set.

Case 2: $w_1\geq 0$, $w_2 \geq 0$, $w_3 = 0$, $z_1=0$, $z_2=0$, $z_3\geq 0$. This leads to
\[
w_1 \leq -6 + \frac{1}{2}z_3 \quad \mbox{and} \quad z_3 \in [2,3].
\]
Hence, $w_1 \geq 0$ cannot be fulfilled and, also in this case, there is no contribution to the solution set.

The Cases 3-6 considering $w_2=0$ all lead to  the contradiction
\[
1 \leq -z_2 \quad \mbox{and} \quad z_2 \geq 0. 
\]
In all these four cases there is no contribution to the solution set.

Case 7: $w_1 = 0$, $w_2 \geq 0$, $w_3 \geq 0$, $z_1\geq 0$, $z_2=0$, $z_3=0$. Then, the inequalities reduce to
\[
\begin{array}{rclcrcl}
 -z_1 & \leq & -6 & \quad \mbox{and} \quad & -6 & \leq & - z_1 \\[1ex]
w_2 -\frac{1}{2}z_1  & \leq & 2 & \quad \mbox{and} \quad & 1 & \leq & w_2 \\[1ex]
w_3 -\frac{1}{2}z_1  & \leq & -2 & \quad \mbox{and} \quad & -3 & \leq & w_3 ,
\end{array}
\]
which leads to
\[
z_1=6, \quad w_2 \in[1,5], \quad w_3 \in [0,1].
\]

So, the contribution to the solution set is the vector
\[
z=\left(
\begin{array}{c}
6  \\[1ex]
0\\[1ex]
0
\end{array}
\right).
\]
\begin{figure}[t]
\unitlength1.0cm
\begin{picture}(12,11)
\includegraphics[width=12cm]{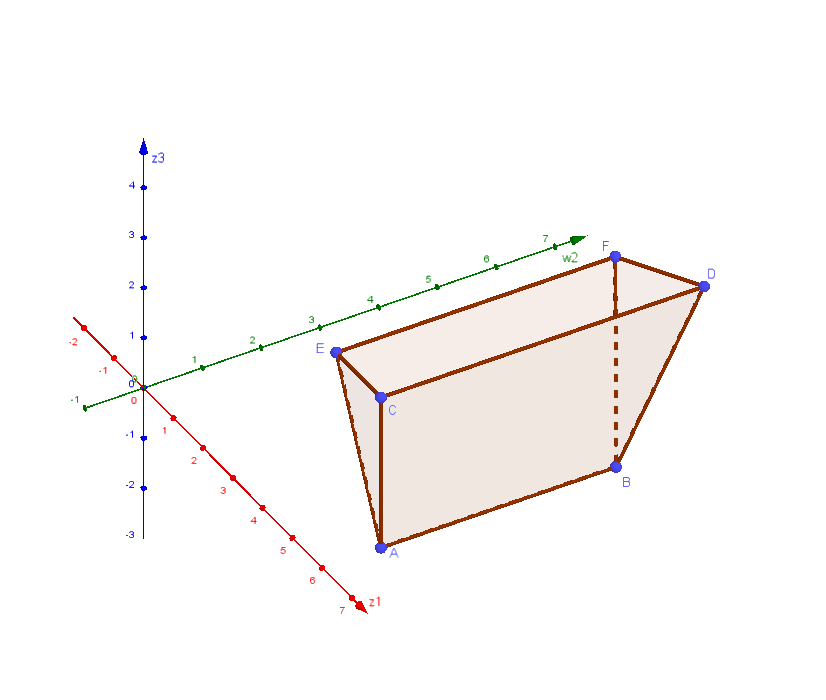}
\end{picture}
\end{figure}

Case 8: $w_1 = 0$, $w_2 \geq 0$, $w_3 = 0$, $z_1\geq 0$, $z_2=0$, $z_3\geq 0$. Then, the inequalities reduce to
\[
\begin{array}{rclcrcl}
-z_1 -\frac{1}{2}z_3 & \leq & -6 & \quad \mbox{and} \quad & -6 & \leq &  - z_1 \\[1ex]
w_2 -\frac{1}{2}z_1  -\frac{1}{2}z_3 & \leq & 2 & \quad \mbox{and} \quad & 1 & \leq & w_2 \\[1ex]
 -\frac{1}{2}z_1  -z_3 & \leq & -2 & \quad \mbox{and} \quad & -3 & \leq &  -z_3
\end{array}
\]
whence
\[
\begin{array}{rclcrcl}
12 & \leq &  2 \cdot z_1 + z_3  & \quad \mbox{and} \quad & z_1 & \leq &  6 \\[1ex]
-4 & \leq & z_1 -2 \cdot w_2 + z_3 &  \quad \mbox{and} \quad & 1 & \leq & w_2 \\[1ex]
 4 & \leq & z_1  + 2 \cdot z_3  & \quad \mbox{and} \quad & z_3 & \leq &  3.
\end{array}
\]
Together with $z_1 \geq 0$, $w_2 \geq 0$, and $z_3 \geq 0$ we get the polyhedron with the vertices
\[
A(6|1|0), \quad B(6|5|0), \quad C(6|1|3), \quad D(6|\frac{13}{2}|3), \quad E(\frac{9}{2}|1|3), \quad F(\frac{9}{2}|\frac{23}{4}|3).
 \]
It is depicted in the figure above.

$A$, $B$, $E$, and $F$ lie in the plane $E_1:\,   2 \cdot z_1 + z_3 =12$ and $B$, $D$, and $F$ lie in the plane
$E_2: \, z_1 -2 \cdot w_2 + z_3 =-4$.

Unifying all eight cases we get the solution set
\[
\Sigma([M],[q])=\left\{
\left(
\begin{array}{c}
z_1 \\[1ex]
z_2 \\[1ex]
z_3
\end{array}
\right) \in  \mathbb{R}^3\, : \, z_2=0, \,  (z_1,z_3) \in L
\right\}
\]
where $L$ is the lightgray shaded region in the following figure.

\psset{xunit=1cm,yunit=1cm, algebraic=true}
\begin{pspicture}(-0.9,-0.9)(8,4)
\psgrid[gridwidth=0.1pt,gridlabels=0,subgriddiv=2,gridcolor=lightgray]
\psaxes[Dx=1,Dy=1,dx=1,dy=1, subticks=2,tickstyle=buttom]{->}(0,0)(-0.9,-0.9)(8,4)
\rput[l](7.7,-0.3){$z_1$}
\rput[l](-0.4,3.8){$z_3$}
\psline*[linewidth=2pt,%
linecolor=lightgray]%
{-}(6,3)(4.5,3)(6,0)
\psplot[linewidth=2pt, linestyle=solid]{4.5}{6}{3}
\psplot[linewidth=2pt, linestyle=solid]{4.5}{6}{-2*x+12}
\psline[linewidth=2pt]{-}%
(6,0)(6,3)
\end{pspicture}

$\underline{M}$ is an M-matrix. Therefore, by Theorem \ref{neuestheorem} $\underline{z}=\sup \Sigma([M],[q])$ where $\underline{z}=\left(
\begin{array}{c}
6  \\[1ex]
0\\[1ex]
3
\end{array}
\right)$ is the unique solution of the LCP defined by $\underline{M}$ and $\underline{q}$. However, $\underline{M}^{-1}>O$ is not fulfilled. Therefore, one cannot expect that
\begin{equation} \label{nichtgutt}
z<\underline{z} \quad \mbox{for any } z \in \Sigma([M],[q]) \backslash \{\underline{z} \}.
\end{equation}
And indeed, (\ref{nichtgutt}) is not true.

Now we consider
the symmetric solution set. The process described in Section 5.2.2 gives  
\[
0\leq z_1^2 -6z_1-z_3^2+3z_3 \quad \mbox{and} \quad 
0 \leq z_3^2 -2 z_3-z_1^2+6z_1
\]
whence
\[
0\leq (z_1-3)^2 -(z_3-\frac{3}{2})^2 - \frac{27}{4} \quad \mbox{and} \quad 
0 \leq (z_3-1)^2 -(z_1-3)^2+8.
\]

So we have to intersect the polyhedron with these hyperbolic cylinders. After projecting the intersection into the $(z_1,z_3)$-plane we get 
the symmetric  solution set
\[
\Sigma_{sym}([M],[q])=\left\{
\left(
\begin{array}{c}
z_1 \\[1ex]
z_2 \\[1ex]
z_3
\end{array}
\right) \in  \mathbb{R}^3\, : \, z_2=0, \,  (z_1,z_3) \in L
\right\}
\]
where $L$ is the lightgray shaded region in the following figure.

\psset{xunit=1cm,yunit=1cm, algebraic=true}
\begin{pspicture}(-0.9,-0.9)(8,4)
\psgrid[gridwidth=0.1pt,gridlabels=0,subgriddiv=2,gridcolor=lightgray]
\psaxes[Dx=1,Dy=1,dx=1,dy=1, subticks=2,tickstyle=buttom]{->}(0,0)(-0.9,-0.9)(8,4)
\rput[l](7.7,-0.3){$z_1$}
\rput[l](-0.4,3.8){$z_3$}
\psline*[linearc=0.1,linewidth=2pt,%
linecolor=lightgray]%
{-}(6,3)(5.62,2)(6,2)
\pscircle*[linecolor=lightgray]%
(5.8,1.9){0.15cm}
\pscircle*[linecolor=lightgray]%
(5.75,1.7){0.15cm}
\pscircle*[linecolor=lightgray]%
(5.7,1.5){0.15cm}
\pscircle*[linecolor=lightgray]%
(5.75,1.3){0.14cm}
\pscircle*[linecolor=lightgray]%
(5.75,1.2){0.1cm}
\pscircle*[linecolor=lightgray]%
(5.75,1.1){0.1cm}
\pscircle*[linecolor=lightgray]%
(5.75,1.0){0.09cm}
\pscircle*[linecolor=lightgray]%
(5.75,0.9){0.09cm}
\pscircle*[linecolor=lightgray]%
(5.77,0.75){0.09cm}
\pscircle*[linecolor=lightgray]%
(5.77,0.7){0.09cm}
\pscircle*[linecolor=lightgray]%
(5.8,0.6){0.09cm}
\psplot[linewidth=2pt, linestyle=solid]{5.83}{7}{1+sqrt((x-3)^2-8)}
\psplot[linewidth=2pt, linestyle=solid]{5.83}{6.5}{1-sqrt((x-3)^2-8)}
\psline[linewidth=2pt]{-}%
(5.83,0.9)(5.83,1.1)
\psplot[linewidth=2pt, linestyle=solid]{-1}{0.17}{1+sqrt((x-3)^2-8)}
\psplot[linewidth=2pt, linestyle=solid]{-0.5}{0.17}{1-sqrt((x-3)^2-8)}
\psline[linewidth=2pt]{-}%
(0.17,0.9)(0.17,1.1)
\psplot[linewidth=2pt, linestyle=solid]{4.5}{6}{3}
\psplot[linewidth=2pt, linestyle=solid]{4.5}{6}{-2*x+12}
\psline[linewidth=2pt]{-}%
(6,0)(6,3)
\psplot[linewidth=2pt, linestyle=solid]{5.6}{6.5}{1.5+sqrt((x-3)^2-27/4)}
\psplot[linewidth=2pt, linestyle=solid]{5.6}{6.5}{1.5-sqrt((x-3)^2-27/4)}
\psline[linewidth=2pt]{-}%
(5.6,1.4)(5.6,1.6)
\psplot[linewidth=2pt, linestyle=solid]{-0.5}{0.4}{1.5+sqrt((x-3)^2-27/4)}
\psplot[linewidth=2pt, linestyle=solid]{-0.5}{0.4}{1.5-sqrt((x-3)^2-27/4)}
\psline[linewidth=2pt]{-}%
(0.4,1.4)(0.4,1.6)
\end{pspicture}

\subsubsection{All $M \in [M]$ are M-matrices}

Let
\[
[M]=\left(
\begin{array}{ccc}
[\frac{1}{3},\frac{1}{2}] & [-\frac{1}{8},-\frac{1}{10}] & [-\frac{1}{8},-\frac{1}{10}] \\[1ex]
[-\frac{1}{8},-\frac{1}{10}] & [\frac{3}{5},\frac{7}{10}] & [-\frac{1}{5},-\frac{1}{6}] \\[1ex]
[-\frac{1}{8},-\frac{1}{10}] & [-\frac{1}{5},-\frac{1}{6}] & [\frac{1}{2},\frac{2}{3}]
\end{array}
\right) \quad \mbox{and} \quad
[q]=\left(
\begin{array}{c}
[-2,4]\\[1ex]
[-2,3] \\[1ex]
[1,2]
\end{array}
\right).
\] 
It is a modified example from \cite{donauer} where the (unsymmetric) solution set was already visua\-lized. All $M \in [M]$ are M-matrices due to 
$\underline{M} 
\left(
\begin{array}{c}
1 \\ 1 \\1
\end{array}
\right)
> o$, Lemma \ref{lemmafan}, and Corollary \ref{korfan}. Furthermore we have
\[
\overline{M}^{-1}=
\frac{1}{911}
\left(
\begin{array}{ccc}
1975 & 375 & 390 \\[1ex]
375 & 1455 & 420 \\[1ex]
390 & 420 & 1530
\end{array}
\right) >O.
\]
So, by Corollary \ref{korfan} we have $ O< \overline{M}^{-1}\leq M^{-1}$ for all $M \in [M]$. The process described in Section 4 gives 
\[
\begin{array}{rclcrcl}
w_1 -\frac{1}{2}z_1+\frac{1}{10}z_2 +\frac{1}{10}z_3 & \leq & 4 & \quad \mbox{and} \quad & -2 & \leq & w_1 -\frac{1}{3}z_1+\frac{1}{8}z_2 +\frac{1}{8}z_3  \\[1ex]
w_2 +\frac{1}{10}z_1-\frac{7}{10}z_2 +\frac{1}{6}z_3 & \leq & 3 & \quad \mbox{and} \quad & -2 & \leq & w_2 +\frac{1}{8}z_1-\frac{3}{5}z_2 +\frac{1}{5}z_3  \\[1ex]
w_3 +\frac{1}{10}z_1+\frac{1}{6}z_2 -\frac{2}{3}z_3 & \leq & 2 & \quad \mbox{and} \quad & 1 & \leq & w_3 +\frac{1}{8}z_1+\frac{1}{5}z_2 -\frac{1}{2}z_3.  
\end{array}
\]

Case 1: $z_1=0$, $z_2=0$, $z_3=0$. We get
\[
w_1 \in [0,4], \quad w_2 \in [0,3], \quad w_3 \in [1,2].
\]
Therefore, the contribution to the solution set is the zero vector.

Case 2: $z_1=0$, $z_2=0$, $w_3=0$. We get 
\[
1\leq -\frac{1}{2}z_3
\]
which is not true for $z_3 \geq 0$. So, there is no contribution to the solution set.

Case 3: $z_1=0$, $w_2=0$, $w_3=0$. We get 
\[
\begin{array}{rclcrcl}
w_1 +\frac{1}{10}z_2 +\frac{1}{10}z_3 & \leq & 4 & \quad \mbox{and} \quad & -2 & \leq & w_1 +\frac{1}{8}z_2 +\frac{1}{8}z_3  \\[1ex]
-\frac{7}{10}z_2 +\frac{1}{6}z_3 & \leq & 3 & \quad \mbox{and} \quad & -2 & \leq & -\frac{3}{5}z_2 +\frac{1}{5}z_3  \\[1ex]
\frac{1}{6}z_2 -\frac{2}{3}z_3 & \leq & 2 & \quad \mbox{and} \quad & 1 & \leq & \frac{1}{5}z_2 -\frac{1}{2}z_3.  
\end{array}
\]
Regarding just the two inequalities
\[
z_3 \geq 3 z_2 -10 , \quad z_3 \leq \frac{2}{5}z_2-2 
\]
we see in the following figure that this cannot be fulfilled for $z_2\geq 0$ and $z_3 \geq 0$.

\psset{xunit=1cm,yunit=1cm, algebraic=true}
\begin{pspicture}(-1.9,-1.9)(8,4)
\psgrid[gridwidth=0.1pt,gridlabels=0,subgriddiv=2,gridcolor=lightgray]
\psaxes[Dx=1,Dy=1,dx=1,dy=1, subticks=2,tickstyle=buttom]{->}(0,0)(-1.9,-1.9)(8,4)
\rput[l](7.7,-0.3){$z_2$}
\rput[l](-0.4,3.8){$z_3$}
\psplot[linewidth=2pt, linestyle=solid]{2.7}{4.5}{3*x-10}
\psplot[linewidth=2pt, linestyle=solid]{1}{8}{0.4*x-2}
\thinlines
\put(4,2){\vector(-3,1){0.5}}
\put(7.5,1){\vector(1,-2){0.3}}
\end{pspicture}

Case 4: $w_1=0$, $z_2=0$, $w_3=0$. We get 
\[
\begin{array}{rclcrcl}
-\frac{1}{2}z_1  +\frac{1}{10}z_3 & \leq & 4 & \quad \mbox{and} \quad & -2 & \leq & -\frac{1}{3}z_1  +\frac{1}{8}z_3  \\[1ex]
w_2 +\frac{1}{10}z_1 +\frac{1}{6}z_3 & \leq & 3 & \quad \mbox{and} \quad & -2 & \leq & w_2  +\frac{1}{8}z_1 +\frac{1}{5}z_3  \\[1ex]
\frac{1}{10}z_1 -\frac{2}{3}z_3 & \leq & 2 & \quad \mbox{and} \quad & 1 & \leq & \frac{1}{8}z_1 -\frac{1}{2}z_3.  
\end{array}
\]
Regarding just the two inequalities
\[
z_3 \geq \frac{8}{3} z_1 -16 , \quad z_3 \leq \frac{1}{4}z_1-2 
\]
we see in the following figure that this cannot be fulfilled for $z_1\geq 0$ and $z_3 \geq 0$.

\psset{xunit=1cm,yunit=1cm, algebraic=true}
\begin{pspicture}(-1.9,-1.9)(9,4)
\psgrid[gridwidth=0.1pt,gridlabels=0,subgriddiv=2,gridcolor=lightgray]
\psaxes[Dx=1,Dy=1,dx=1,dy=1, subticks=2,tickstyle=buttom]{->}(0,0)(-1.9,-1.9)(9,4)
\rput[l](8.7,-0.3){$z_1$}
\rput[l](-0.4,3.8){$z_3$}
\psplot[linewidth=2pt, linestyle=solid]{5.5}{7.5}{2.66*x-16}
\psplot[linewidth=2pt, linestyle=solid]{1}{9}{0.25*x-2}
\thinlines
\put(6,0){\vector(-3,1){0.5}}
\put(8,0){\vector(1,-2){0.3}}
\end{pspicture}

Case 5: $z_1=0$, $w_2=0$, $z_3=0$. We get 
\[
\begin{array}{rclcrcl}
w_1 +\frac{1}{10}z_2  & \leq & 4 & \quad \mbox{and} \quad & -2 & \leq & w_1 +\frac{1}{8}z_2  \\[1ex]
-\frac{7}{10}z_2  & \leq & 3 & \quad \mbox{and} \quad & -2 & \leq & -\frac{3}{5}z_2   \\[1ex]
w_3 +\frac{1}{6}z_2  & \leq & 2 & \quad \mbox{and} \quad & 1 & \leq & w_3+ \frac{1}{5}z_2 .  
\end{array}
\]
Regarding  the two inequalities
\[
w_1 \leq   -\frac{1}{10} z_2 +4 , \quad w_1 \geq -\frac{1}{8}z_2-2 \quad \mbox{for } z_2 \in [0,\frac{10}{3}] 
\]
we get for $w_1\geq 0$ and $z_2 \geq 0$

\psset{xunit=1cm,yunit=1cm, algebraic=true}
\begin{pspicture}(-0.9,-0.9)(7,5)
\psgrid[gridwidth=0.1pt,gridlabels=0,subgriddiv=2,gridcolor=lightgray]
\psline*[linewidth=2pt,%
linecolor=lightgray]%
{-}(0,0)(3.33,0)(3.33,3.66)(0,4)
\psaxes[Dx=1,Dy=1,dx=1,dy=1, subticks=2,tickstyle=buttom]{->}(0,0)(-0.9,-0.9)(7,5)
\rput[l](6.7,-0.3){$z_2$}
\rput[l](-0.5,4.6){$w_1$}
\psplot[linewidth=2pt, linestyle=solid]{0}{3.33}{-0.1*x+4}
\psline[linewidth=2pt]{-}%
(3.33,0)(3.33,3.66)
\end{pspicture}

Regarding  the two inequalities
\[
w_3 \leq   -\frac{1}{6} z_2 +2 , \quad w_3 \geq -\frac{1}{5}z_2+1 \quad \mbox{for } z_2 \in [0,\frac{10}{3}] 
\]
we get for $w_3\geq 0$ and $z_2 \geq 0$

\psset{xunit=1cm,yunit=1cm, algebraic=true}
\begin{pspicture}(-0.9,-0.9)(7,3)
\psgrid[gridwidth=0.1pt,gridlabels=0,subgriddiv=2,gridcolor=lightgray]
\psline*[linewidth=2pt,%
linecolor=lightgray]%
{-}(0,1)(3.33,0.33)(3.33,1.44)(0,2)
\psaxes[Dx=1,Dy=1,dx=1,dy=1, subticks=2,tickstyle=buttom]{->}(0,0)(-0.9,-0.9)(7,3)
\rput[l](6.7,-0.3){$z_2$}
\rput[l](-0.5,2.6){$w_3$}
\psplot[linewidth=2pt, linestyle=solid]{0}{3.33}{-0.16*x+2}
\psplot[linewidth=2pt, linestyle=solid]{0}{3.33}{-0.2*x+1}
\psline[linewidth=2pt]{-}%
(3.33,1.44)(3.33,0.33)
\end{pspicture}

So the contribution to the solution set is
\[
\left\{
\left(
\begin{array}{c}
z_1 \\[1ex]
z_2 \\[1ex]
z_3
\end{array}
\right) \in  \mathbb{R}^3\, : \, z_1=0, \, z_2 \in [0,\frac{10}{3}], \,  z_3=0 
\right\}.
\]

Case 6: $w_1=0$, $z_2=0$, $z_3=0$. We get 
\[
\begin{array}{rclcrcl}
-\frac{1}{2}z_1  & \leq & 4 & \quad \mbox{and} \quad & -2 & \leq & -\frac{1}{3}z_1  \\[1ex]
w_2+\frac{1}{10}z_1  & \leq & 3 & \quad \mbox{and} \quad & -2 & \leq &w_2 +\frac{1}{8}z_1   \\[1ex]
w_3 +\frac{1}{10}z_1  & \leq & 2 & \quad \mbox{and} \quad & 1 & \leq & w_3+ \frac{1}{8}z_1 .  
\end{array}
\]
Regarding  the two inequalities
\[
w_2 \leq   -\frac{1}{10} z_1 +3 , \quad w_2 \geq -\frac{1}{8}z_1-2 \quad \mbox{for } z_1 \in [0,6] 
\]
we get for $z_1\geq 0$ and $w_2 \geq 0$

\psset{xunit=1cm,yunit=1cm, algebraic=true}
\begin{pspicture}(-0.9,-0.9)(7,5)
\psgrid[gridwidth=0.1pt,gridlabels=0,subgriddiv=2,gridcolor=lightgray]
\psline*[linewidth=2pt,%
linecolor=lightgray]%
{-}(0,0)(6,0)(6,2.4)(0,3)
\psaxes[Dx=1,Dy=1,dx=1,dy=1, subticks=2,tickstyle=buttom]{->}(0,0)(-0.9,-0.9)(7,5)
\rput[l](6.7,-0.3){$z_1$}
\rput[l](-0.5,4.6){$w_2$}
\psplot[linewidth=2pt, linestyle=solid]{0}{6}{-0.1*x+3}
\psline[linewidth=2pt]{-}%
(6,0)(6,2.4)
\end{pspicture}

Regarding  the two inequalities
\[
w_3 \leq   -\frac{1}{10} z_1 +2 , \quad w_3 \geq -\frac{1}{8}z_1+1 \quad \mbox{for } z_1 \in [0,6] 
\]
we get for $z_1\geq 0$ and $w_3\geq 0$

\psset{xunit=1cm,yunit=1cm, algebraic=true}
\begin{pspicture}(-0.9,-0.9)(7,3)
\psgrid[gridwidth=0.1pt,gridlabels=0,subgriddiv=2,gridcolor=lightgray]
\psline*[linewidth=2pt,%
linecolor=lightgray]%
{-}(0,1)(6,0.25)(6,1.4)(0,2)
\psaxes[Dx=1,Dy=1,dx=1,dy=1, subticks=2,tickstyle=buttom]{->}(0,0)(-0.9,-0.9)(7,3)
\rput[l](6.7,-0.3){$z_1$}
\rput[l](-0.5,2.6){$w_3$}
\psplot[linewidth=2pt, linestyle=solid]{0}{6}{-0.1*x+2}
\psplot[linewidth=2pt, linestyle=solid]{0}{6}{-x/8+1}
\psline[linewidth=2pt]{-}%
(6,0.25)(6,1.4)
\end{pspicture}

So the contribution to the solution set is
\[
\left\{
\left(
\begin{array}{c}
z_1 \\[1ex]
z_2 \\[1ex]
z_3
\end{array}
\right) \in  \mathbb{R}^3\, :  \, z_1 \in [0,6], \, z_2=0,\,   z_3=0 
\right\}.
\]

Case 7: $w_1 = 0$, $w_2 =0$, $z_3 = 0$. Then, the inequalities reduce to
\[
\begin{array}{rclcrcl}
-\frac{1}{2}z_1 +\frac{1}{10}z_2 & \leq & 4 & \quad \mbox{and} \quad & -2 & \leq &  -\frac{1}{3} z_1+\frac{1}{8}z_2 \\[1ex]
\frac{1}{10}z_1  -\frac{7}{10}z_2 & \leq & 3 & \quad \mbox{and} \quad & -2 & \leq & \frac{1}{8}z_1-\frac{3}{5}z_2 \\[1ex]
w_3+ \frac{1}{10}z_1  +\frac{1}{6}z_2 & \leq & 2 & \quad \mbox{and} \quad & 1 & \leq & w_3 +\frac{1}{8}z_1+\frac{1}{5}z_2
\end{array}
\]
whence
\[
\begin{array}{rclcrcl}
z_2 & \leq &  5 \cdot z_1 + 40  & \quad \mbox{and} \quad & z_2 & \geq &  \frac{8}{3}z_1-16 \\[1ex]
z_2 & \geq & \frac{1}{7}z_1 - \frac{30}{7} &  \quad \mbox{and} \quad & z_2 & \leq & \frac{5}{24}z_1+\frac{10}{3} \\[1ex]
60 & \geq & 30 \cdot w_3 + 3 z_1 + 5 z_2   & \quad \mbox{and} \quad & 40 & \leq &  40 \cdot w_3+5z_1+8z_2.
\end{array}
\]
Regarding first only the first four inequalities we get 

\psset{xunit=1cm,yunit=1cm, algebraic=true}
\begin{pspicture}(-0.9,-0.9)(9,6)
\psgrid[gridwidth=0.1pt,gridlabels=0,subgriddiv=2,gridcolor=lightgray]
\psline*[linewidth=2pt,%
linecolor=lightgray]%
{-}(0,0)(6,0)(7.86,4.97)(0,3.33)
\psaxes[Dx=1,Dy=1,dx=1,dy=1, subticks=2,tickstyle=buttom]{->}(0,0)(-0.9,-0.9)(9,6)
\rput[l](8.7,-0.3){$z_1$}
\rput[l](-0.5,5.6){$z_2$}
\psplot[linewidth=2pt, linestyle=solid]{6}{7.86}{2.66*x-16}
\psplot[linewidth=2pt, linestyle=solid]{0}{7.87}{0.21*x+3.33}
\put(8,5){$V$}
\put(7.85,4.95){\circle{0.05}}
\put(7.85,4.95){\circle{0.1}}
\put(7.85,4.95){\circle{0.075}}
\put(7.85,4.95){\circle{0.02}}

\end{pspicture}

taking into account  $z_1 \geq 0$ and $z_2 \geq 0$ with $V(\frac{464}{59}|\frac{880}{177})$. 

Together with the last two inequalities we get the polyhedron with the vertices
\[
V_1(0|0|2), \quad V_2(0|0|1), \quad V_3(0|\frac{10}{3}|\frac{1}{3}), \quad V_4(0|\frac{10}{3}|\frac{13}{9}), \quad V_5(6|0|\frac{42}{30}), \quad V_6(6|0|\frac{1}{4}),
\]
\[
 V_7(\frac{464}{59}|\frac{880}{177}|0), \quad  V_8(\frac{464}{59}|\frac{880}{177}|\frac{1022}{2655}), \quad V_9(\frac{504}{79}|\frac{80}{79}|0), \quad
V_{10}(2|\frac{15}{4}|0) .
 \]
It is depicted in the following figure.



\begin{figure}[h]
\unitlength1.0cm
\begin{picture}(12,14)
\includegraphics[scale=1.3]{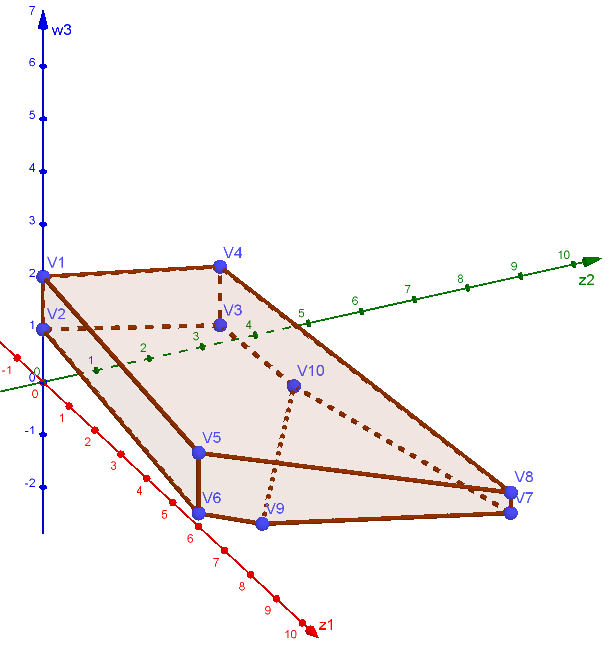}
\end{picture}
\end{figure}

$V_1$, $V_4$, $V_5$, and $V_8$ lie in the plane $E_1:\,   3\cdot z_1 +5 \cdot z_2+30 \cdot w_3=60$ and $V_2$, $V_3$, $V_6$, $V_9$, and $V_{10}$ lie in the plane
$E_2: \,  5 z_1 +8 z_2+40 \cdot w_3=40$.

So the contribution to the solution set is
\[
\left\{
\left(
\begin{array}{c}
z_1 \\[1ex]
z_2 \\[1ex]
z_3
\end{array}
\right) \in  \mathbb{R}^3\, :  \, (z_1, z_2) \in L ,\,   z_3=0 
\right\}
\]
where $L$ is the figure one before.

Now we consider the symmetric solution set. The process described in Section 5.2.2 gives
\[
0\leq \frac{7}{10}z_2^2 +3z_2-\frac{1}{3}z_1^2+2z_1 \quad \mbox{and} \quad 
0 \leq \frac{1}{2}z_1^2+4z_1 -\frac{3}{5} z_2^2+2z_2
\]
which leads to
\[
0\leq \frac{7}{10}(z_2+\frac{15}{7})^2 -\frac{1}{3}(z_1-3)^2 - \frac{3}{14} \quad \mbox{and} \quad 
0 \leq \frac{1}{2}(z_1+4)^2 -\frac{3}{5}(z_2-\frac{5}{3})^2-\frac{19}{3}.
\]
we have to intersect the polyhedron above with the corresponding  hyperbolic cylinders.
Using GeoGebra, for example, we see that the solution set is not truncated by the hyperbola.\newline
\includegraphics[scale=0.53]{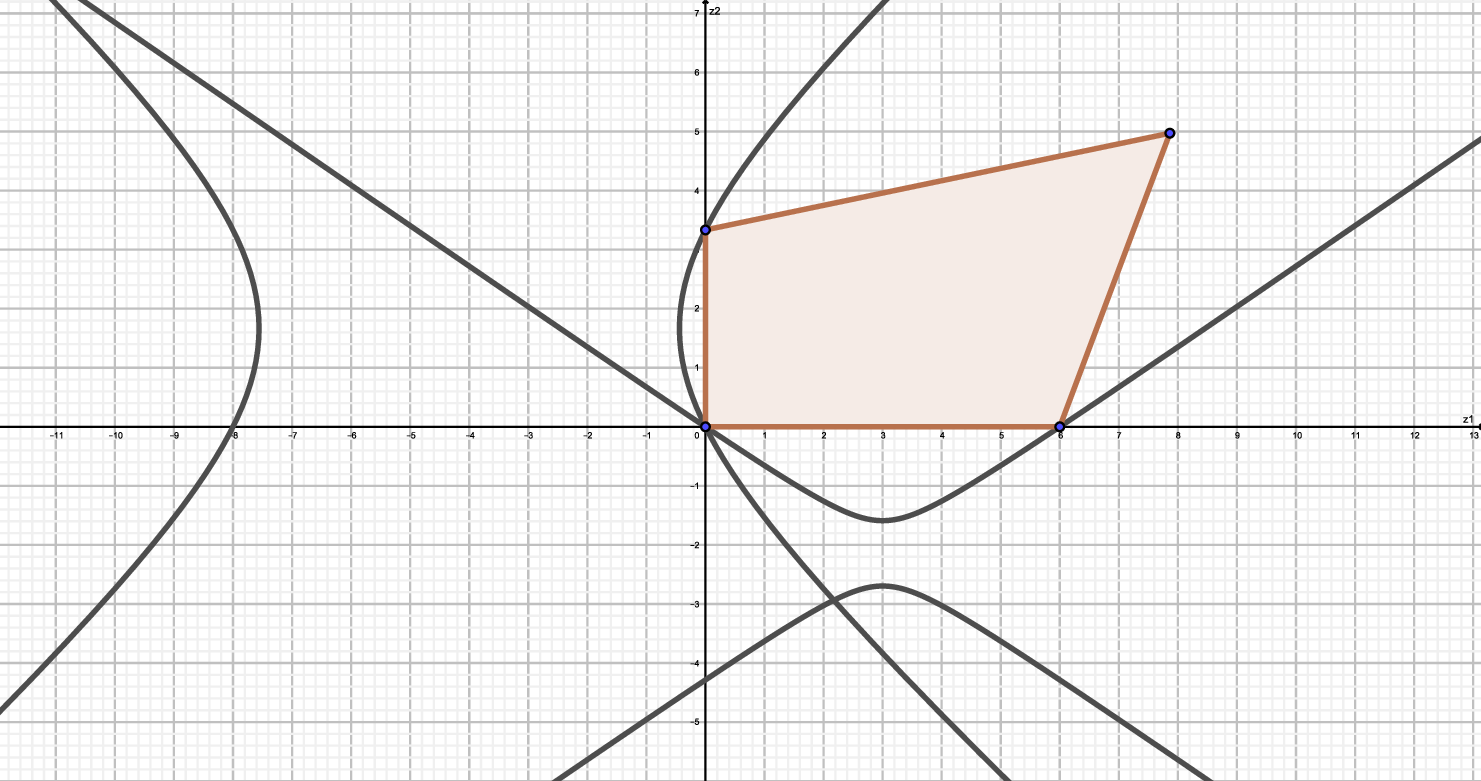}

Case 8: $w_1 = 0$, $w_2 =0$, $w_3 = 0$. Then, the inequalities reduce to
\[
\begin{array}{rclcrcl}
-\frac{1}{2}z_1 +\frac{1}{10}z_2 + \frac{1}{10}z_3 & \leq & 4 & \quad \mbox{and} \quad & -2 & \leq &  -\frac{1}{3} z_1+\frac{1}{8}z_2+\frac{1}{8}z_3 \\[1ex]
\frac{1}{10}z_1  -\frac{7}{10}z_2  +\frac{1}{6}z_3& \leq & 3 & \quad \mbox{and} \quad & -2 & \leq & \frac{1}{8}z_1-\frac{3}{5}z_2 +\frac{1}{5}z_3  \\[1ex]
\frac{1}{10}z_1  +\frac{1}{6}z_2  -\frac{2}{3}z_3  & \leq & 2 & \quad \mbox{and} \quad & 1 & \leq & \frac{1}{8}z_1+\frac{1}{5}z_2  -\frac{1}{2}z_3 
\end{array}
\]
whence
\[
\begin{array}{rclcrcl}
-5z_1 +z_2 + z_3 & \leq & 40 & \quad \mbox{and} \quad & -48 & \leq &  -8z_1+3z_2+3z_3 \\[1ex]
3z_1  -21 z_2  +5 z_3& \leq & 90 & \quad \mbox{and} \quad & -80 & \leq & 5z_1-24z_2 +8z_3  \\[1ex]
3z_1  +5z_2  - 20 z_3  & \leq & 60 & \quad \mbox{and} \quad & 40 & \leq & 5z_1+8z_2  -20z_3. 
\end{array}
\]
The contribution to the solution set is the pyramid with apex 
\[
S(\frac{5712}{607}/\frac{3790}{607}/\frac{1730}{607})
\]
and the triangle ABV
with
\[
A(\frac{504}{79}/\frac{80}{79}/0), \quad B(2/\frac{15}{4}/0), \quad V(\frac{464}{59}/\frac{880}{177}/0)
\]
as its base.  See the figure above.

\begin{figure}[h]
\unitlength1.0cm
\begin{picture}(12,10.5)
\includegraphics[scale=0.7]{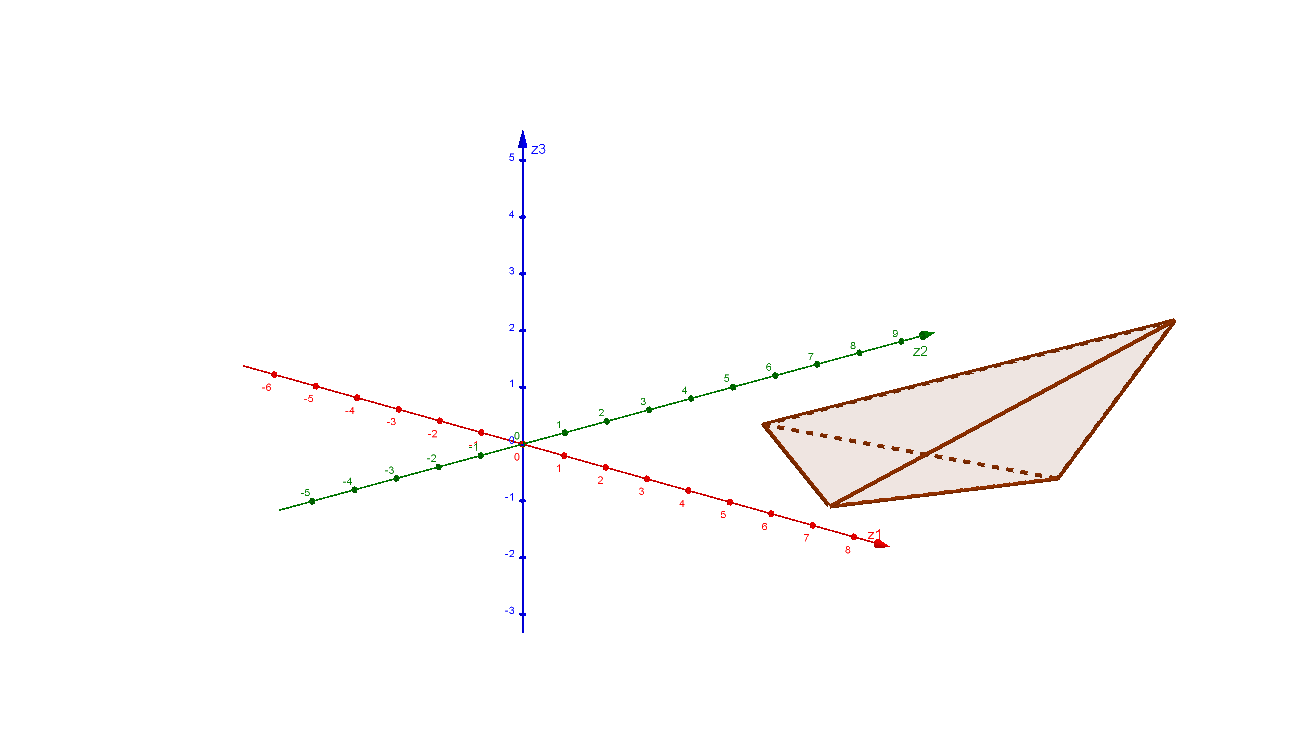}
\end{picture}
\end{figure}

Unifying all 8 cases we get as the solution set the following figure.

\begin{figure}[h]
\unitlength1.0cm
\begin{picture}(12,10.5)
\includegraphics[scale=0.7]{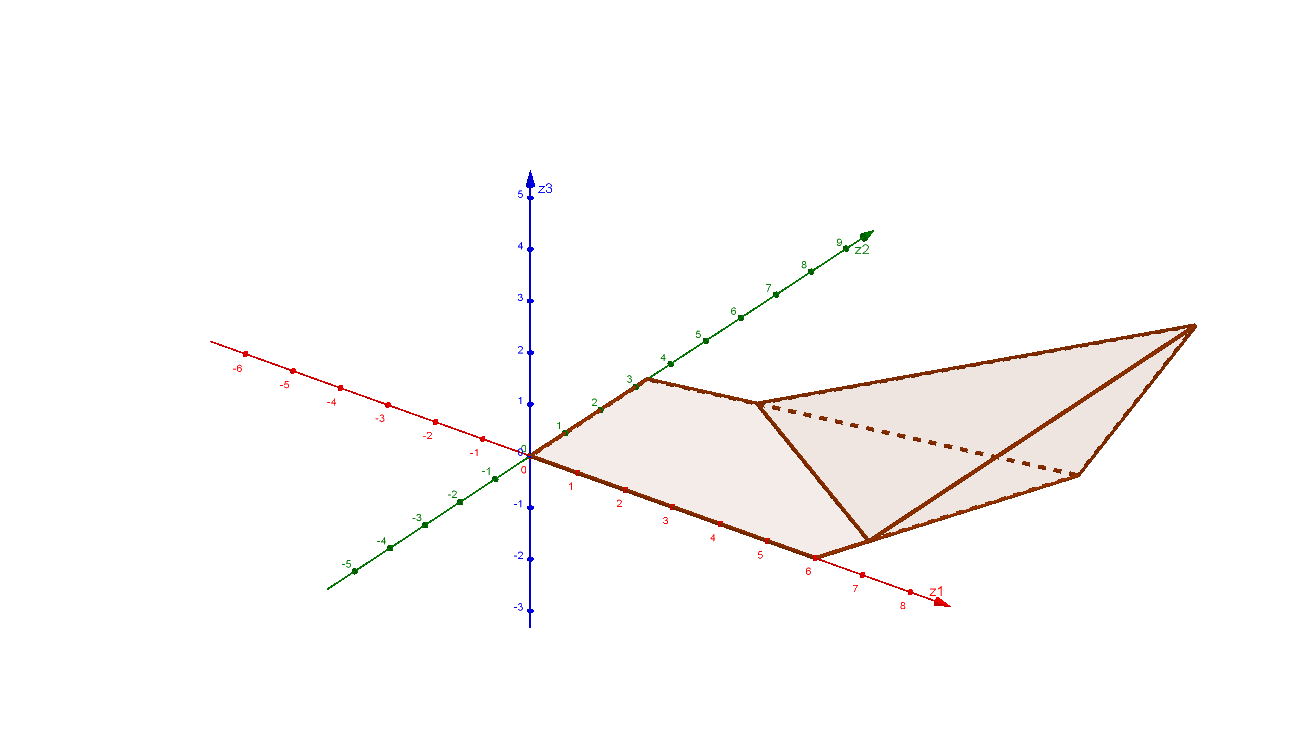}
\end{picture}
\end{figure}

We will show, that in this example the solution set and the symmetric solution set are identical. To gain this we have to show that the pyramid from Case 8 is not truncated by the quadrics from the Fourier-Motzkin elimination process. We consider the three planes
\[
\begin{array}{c}
E_2: \quad \underline{q}_1 = -\underline{m}_{11} z_1 -\underline{m}_{12} z_2 -\underline{m}_{13}z_3 \\[1ex]
E_4: \quad \underline{q}_2 = -\underline{m}_{21} z_1 -\underline{m}_{22} z_2 -\underline{m}_{23}z_3 \\[1ex]
E_6: \quad \underline{q}_3 = -\underline{m}_{31} z_1 -\underline{m}_{32} z_2 -\underline{m}_{33}z_3 
\end{array}
\] 
to which the lateral faces of the pyramid belong. (Note, that the vector
\[
\underline{z}=
\left(
\begin{array}{c}
\displaystyle{\frac{5712}{607}}\\[2ex]
\displaystyle{\frac{3790}{607}}\\[2ex]
\displaystyle{\frac{1730}{607}}
\end{array}
\right)
\]
corresponding to the apex is the unique solution to the LCP defined by $\underline{M}$ and $\underline{q}$. See Theorem \ref{neuestheorem}.) 

The next aim will be to show, that the edge from $A$ to $S$ as well as the edge from $B$ to $S$ as well as the edge from $V$ to $S$ all belong to the symmetric solution set.

Firstly, the vector 
\[
z=
\left(
\begin{array}{c}
2\\[1ex]
\displaystyle{\frac{15}{4}}\\[1ex]
0
\end{array}
\right) \quad \mbox{ (corresponding to the point $B$)}
\]
is the unique solution of the LCP defined by the symmetric matrix
\[
M=
\left(
\begin{array}{rrr}
\displaystyle{\frac{1}{2}} & -\displaystyle{\frac{1}{8}} & - \displaystyle{\frac{1}{8}} \\[2ex]
-\displaystyle{\frac{1}{8}} & \displaystyle{\frac{3}{5}} & - \displaystyle{\frac{1}{5}} \\[2ex]
-\displaystyle{\frac{1}{8}} & -\displaystyle{\frac{1}{5}} &  \displaystyle{\frac{1}{2}} 
\end{array}
\right)
\quad \mbox{and} \quad
q=\left(
\begin{array}{c}
 -\displaystyle{ \frac{17}{32}} \\[2ex]
-2 \\[1ex]
1 
\end{array}
\right) \mbox{ with } q=-Mz.
\]
Defining the symmetric matrix
\[
M(t)=M-
\left(
\begin{array}{rrr}
t &0 & 0 \\[1ex]
0 & 0 & 0 \\[1ex]
0 & 0&  0 
\end{array}
\right)
\quad \mbox{with $t\in (0,\frac{1}{6})$} \quad \mbox{and} \quad
q(s)=q-\left(
\begin{array}{c}
 s \\[1ex]
0 \\[1ex]
0 
\end{array}
\right)  \quad \mbox{with $s\in (0,\frac{47}{32})$ }
\]
we have $\underline{M} \leq M(t) \leq M$ and $\underline{q} \leq q(s)\leq q$. Let $z(t,s)$ be the unique solution of the LCP defined by  $M(t)$ and $q(s)$. Then we have
\[
o \leq z < z(t,s) = -M(t)^{-1}\cdot q(s) < \underline{z}
\]
by (the second part of) Theorem \ref{neuestheorem}. Therefore, $z(t,s)\in E_4 \cap E_6$ and
$z(t,s)$ runs from $z$ to $\underline{z}$ as soon as $t \to \frac{1}{6}$ and $s \to \frac{47}{32}$; i.e., the edge from B to S belongs to the symmetric solution set. 

Secondly, the vector 
\[
z=
\left(
\begin{array}{c}
\displaystyle{\frac{504}{79}}\\[2ex]
\displaystyle{\frac{80}{79}}\\[2ex]
0
\end{array}
\right)  \quad \mbox{ (corresponding to the point $A$)}
\]
is the unique solution of the LCP defined by the symmetric matrix
\[
M=
\left(
\begin{array}{rrr}
\displaystyle{\frac{1}{3}} & -\displaystyle{\frac{1}{8}} & - \displaystyle{\frac{1}{8}} \\[2ex]
-\displaystyle{\frac{1}{8}} & \displaystyle{\frac{7}{10}} & - \displaystyle{\frac{1}{5}} \\[2ex]
-\displaystyle{\frac{1}{8}} & -\displaystyle{\frac{1}{5}} &  \displaystyle{\frac{1}{2}} 
\end{array}
\right)
\quad \mbox{and} \quad
q=\left(
\begin{array}{c}
-2\\[2ex]
 \displaystyle{ \frac{7}{79}} \\[2ex]
1 
\end{array}
\right)   \quad \mbox{with } q=-Mz  .
\]
Defining the symmetric matrix
\[
M(t)=M-
\left(
\begin{array}{rrr}
0 &0 & 0 \\[1ex]
0 & t & 0 \\[1ex]
0 & 0&  0 
\end{array}
\right)
\quad \mbox{with $t\in (0,\frac{1}{10})$} \quad \mbox{and} \quad
q(s)=q-\left(
\begin{array}{c}
 0 \\[1ex]
s \\[1ex]
0 
\end{array}
\right)
\quad \mbox{with $s\in (0,\frac{165}{79})$ }
\]
we have $\underline{M} \leq M(t) \leq M$ and $\underline{q} \leq q(s)\leq q$. Let $z(t,s)$ be the unique solution of the LCP defined by  $M(t)$ and $q(s)$.  Then we have
\[
o \leq z < z(t,s) = -M(t)^{-1}\cdot q(s)<\underline{z}
\]
by (the second part of) Theorem \ref{neuestheorem}. Therefore, $z(t,s)\in E_2 \cap E_6$ and
$z(t,s)$ runs from $z$ to $\underline{z}$ as soon as $t \to \frac{1}{10}$ and $s \to \frac{165}{79}$; i.e., the edge from A to S belongs to the symmetric solution set. 

Thirdly, the vector 
\[
z=
\left(
\begin{array}{c}
\displaystyle{\frac{464}{59}}\\[2ex]
\displaystyle{\frac{880}{177}}\\[2ex]
0
\end{array}
\right)  \quad \mbox{ (corresponding to the point $V$)}
\]
is the unique solution of the LCP defined by the symmetric matrix $\underline{M}$ and
\[
q=\left(
\begin{array}{c}
-2\\[1ex]
-2\\[1ex]
 \displaystyle{ \frac{350}{177}} 
\end{array}
\right)   \quad \mbox{with } q=-\underline{M}z.
\]
Defining 
\[
q(s)=q-\left(
\begin{array}{c}
 0 \\[1ex]
0 \\[1ex]
s 
\end{array}
\right)  \mbox{with $s\in (0,\frac{173}{177})$ }
\]
we have $\underline{q} \leq q(s)\leq q$. Let $z(s)$ be the unique solution of the LCP defined by  $\underline{M}$ and $q(s)$.  Then we have
\[
o \leq z < z(s) = -\underline{M}^{-1}\cdot q(s) < \underline{z}
\]
by (the second part of) Theorem \ref{neuestheorem}. Therefore, $z(s)\in E_2 \cap E_4$ and
$z(s)$ runs from $z$ to $\underline{z}$ as soon as $s \to \frac{173}{177}$; i.e., the edge from V to S belongs to the symmetric solution set. 

In Case 7 we have already shown, that the base of the pyramid belongs to the symmetric solution set. So, if a quadric from the Fourier-Motzkin elimination process truncates a lateral face from the pyramid without truncating the edges and the base, the resulting curve from the intersection of a triangle and a quadric (without truncating the edges of the triangle) must be an ellipse  (or a circle). Therefore, we have the situation in the following figure.

\includegraphics[scale=0.7]{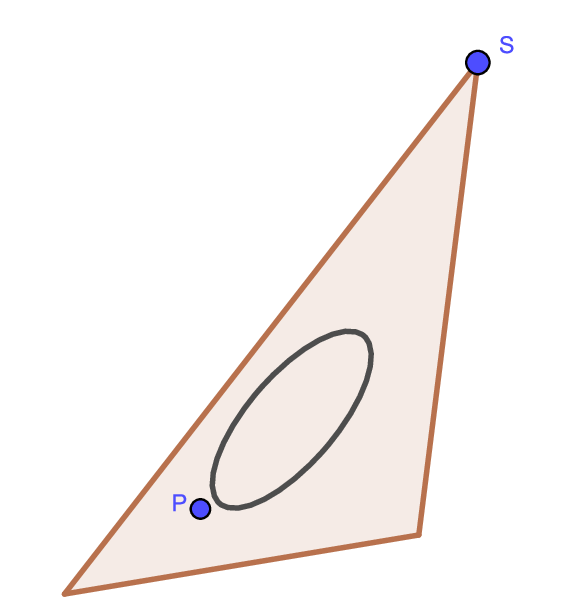}

The corresponding vector, say $z_P$, to the point $P$ outside the ellipse belongs to the symmetric solution set, but all the vectors that correspond to the points within the ellipse do not. That means, that there exist some symmetric matrix $M \in [M]$ and some $q \in [q]$ such that $z_P>o $ solves the corresponding LCP. W.l.o.g. let the considered lateral face be part of the plane
\[
E_2: \quad \underline{q}_1=-\underline{m}_{11} \cdot z_1
- \underline{m}_{12} \cdot z_2 - \underline{m}_{13}\cdot z_3.
\]
Then
\[
M=\left(
\begin{array}{rrr}
\underline{m}_{11} & \underline{m}_{12} & \underline{m}_{13} \\[1ex]
\underline{m}_{12} & m_{22} & m_{23} \\[1ex]
\underline{m}_{13} & m_{23} & m_{33} 
\end{array}
\right)
\in [M] \quad \mbox{and} \quad
q=\left(
\begin{array}{c}
 \underline{q}_{1} \\[1ex]
q_2 \\[1ex]
q_3
\end{array}
\right)\in [q].
\]
Defining the symmetric matrix
\[
M(t_1,t_2,t_3)=M-
\left(
\begin{array}{rrr}
0 &0 & 0 \\[1ex]
0 & t_1 & t_2 \\[1ex]
0 & t_2 &  t_3 
\end{array}
\right)
\quad \mbox{with} \quad
\left\{
\begin{array}{r}
t_1 \in (0,m_{22}-\underline{m}_{22})\\[1ex]
t_2 \in (0,m_{23}-\underline{m}_{23})\\[1ex]
t_3 \in (0,m_{33}-\underline{m}_{33}) 
\end{array}
\right.
\]
and
\[
q(s_1,s_2)=q-\left(
\begin{array}{c}
 0 \\[1ex]
s_1 \\[1ex]
s_2
\end{array}
\right)
\quad \mbox{with $s_1\in (0,q_{2}-\underline{q}_{2})$, \, $s_2\in (0,q_{3}-\underline{q}_{3})$  }
\]
we have $\underline{M} \leq M(t_1,t_2,t_3) \leq M$ and $\underline{q} \leq q(s_1,s_2)\leq q$. Let $z(t_1,t_2,t_3,s_1,s_2)$ be the unique solution of the LCP defined by  $M(t_1,t_2,t_3)$ and $q(s_1,s_2)$.  Then we have
\begin{equation} \label{bleistiftstern}
o < z_P < z(t_1,t_2,t_3,s_1,s_2)    
= -M(t_1,t_2,t_3)^{-1}\cdot q(s_1,s_2)<\underline{z}
\end{equation}
by (the second part of) Theorem \ref{neuestheorem}, and all the vectors $z_P$, $z(t_1,t_2,t_3,s_1,s_2)$, $\underline{z}$ belong to $\Sigma_{sym}([M],[q])$. This means, that
the path from $P$ to $S$ via $ z(t_1,t_2,t_3,s_1,s_2)$ does not cross the ellipse, which is not possible due to (\ref{bleistiftstern}) and the figure above.

\subsubsection{All $M \in [M]$ are $H_+$-matrices}

Let
\[
[M]=\left(
\begin{array}{ccc}
[4,5] & [-1,2] & 0 \\[1ex]
[-1,2] & [2,3] & 1 \\[1ex]
0 & 1 & 2
\end{array}
\right) \quad \mbox{and} \quad
[q]=\left(
\begin{array}{c}
[-2,-1] \\[1ex]
[-1,1] \\[1ex]
[-2,-1]
\end{array}
\right).
\] 
Since
\[
\left(
\begin{array}{rrr}
4 & -2 & 0 \\[1ex]
-2 & 2 & -1 \\[1ex]
0 & -1 & 2
\end{array}
\right)^{-1}=
\left(
\begin{array}{rrr}
\displaystyle{\frac{3}{4}} & 1 & \displaystyle{\frac{1}{2}} \\[1ex]
1 & 2 & 1 \\[1ex]
\displaystyle{\frac{1}{2}} & 1 & 1
\end{array}
\right)\geq O
\]
we get via Lemma \ref{lemmafan} and Corollary \ref{korfan} that all $M\in [M]$ are $H_{+}$-matrices. The process described in Section 4 gives
\[
\begin{array}{rclcrcl}
w_1 -5z_1-2z_2  & \leq & -1 & \quad \mbox{and} \quad & -2 & \leq & w_1 -4z_1+z_2   \\[1ex]
w_2 -2z_1-3z_2 -z_3 & \leq & 1 & \quad \mbox{and} \quad & -1 & \leq & w_2 +z_1-2z_2 -z_3  \\[1ex]
w_3 -z_2 -2\cdot z_3 & \leq & -1 & \quad \mbox{and} \quad & -2 & \leq & w_3 -z_2 -2\cdot z_3.  
\end{array}
\]
We have to consider eight cases.\\
Case 1: $z_1=0$, $z_2=0$, $z_3=0$,  $w_1\geq 0$, $w_2\geq 0$, $w_3\geq 0$. \\
Case 2: $z_1=0$, $w_2=0$, $z_3=0$, $w_1\geq 0$, $z_2\geq 0$, $w_3\geq 0$.\\
Case 3: $z_1=0$, $z_2=0$, $w_3=0$, $w_1\geq 0$, $w_2\geq 0$, $z_3\geq 0$.\\
Case 4: $z_1=0$, $w_2=0$, $w_3=0$, $w_1\geq 0$, $z_2\geq 0$, $z_3\geq 0$. \\
Case 5: $w_1=0$, $z_2=0$, $z_3=0$, $z_1\geq 0$, $w_2\geq 0$, $w_3\geq 0$.\\
Case 6: $w_1=0$, $w_2=0$, $z_3=0$, $z_1\geq 0$, $z_2\geq 0$, $w_3\geq 0$. \\
Case 7: $w_1=0$, $z_2=0$, $w_3=0$, $z_1\geq 0$, $w_2\geq 0$, $z_3\geq 0$. \\
Case 8: $w_1=0$, $w_2=0$, $w_3=0$, $z_1\geq 0$, $z_2\geq 0$, $z_3\geq 0$. \\
Case 1 and Case 3 give no contribution to the solution set, since the first inequality reduces to $w_1 \leq -1$. Also Case 5 ends with a contradiction, namely $w_3 \leq -1$.

Case 2: $z_1=0$, $w_2=0$, $z_3=0$, $w_1\geq 0$, $z_2\geq 0$, $w_3\geq 0$. Then the inequalities reduce to
\[
\begin{array}{rclcrcl}
w_1 -2z_2  & \leq & -1 & \quad \mbox{and} \quad & -2 & \leq & w_1 +z_2   \\[1ex]
-3z_2  & \leq & 1 & \quad \mbox{and} \quad & -1 & \leq & -2z_2  \\[1ex]
w_3 -z_2  & \leq & -1 & \quad \mbox{and} \quad & -2 & \leq & w_3 -z_2.  
\end{array}
\]
In particular, we get $z_2 \in [0,\frac{1}{2}]$ and $ w_3 +1 \leq z_2$, which is a contradiction. So, Case 2 gives no contribution to the solution set.

Case 4: $z_1=0$, $w_2=0$, $w_3=0$, $w_1\geq 0$, $z_2\geq 0$, $z_3\geq 0$. Then the inequalities reduce to
\[
\begin{array}{rclcrcl}
w_1 -2z_2  & \leq & -1 & \quad \mbox{and} \quad & -2 & \leq & w_1 +z_2   \\[1ex]
-3z_2 -z_3 & \leq & 1 & \quad \mbox{and} \quad & -1 & \leq & -2z_2 -z_3  \\[1ex]
-z_2 -2\cdot z_3 & \leq & -1 & \quad \mbox{and} \quad & -2 & \leq &  -z_2 -2\cdot z_3,  
\end{array}
\]
whence
\[
\begin{array}{rclcrcl}
w_1  & \leq & 2z_2 -1 & \quad \mbox{and} \quad & -z_2 -2 & \leq & w_1    \\[1ex]
-3z_2-1  & \leq & z_3 & \quad \mbox{and} \quad & z_3 & \leq & -2z_2 +1  \\[1ex]
 -\frac{1}{2}z_2 +\frac{1}{2} & \leq & z_3 & \quad \mbox{and} \quad & z_3 & \leq & -\frac{1}{2}z_2 +1.  
\end{array}
\]
Regarding  the inequalities
\[
 z_3  \leq  -2z_2 +1 , \quad  -\frac{1}{2}z_2 +\frac{1}{2} \leq  z_3 , \quad z_2\geq 0, \quad z_3 \geq 0
\]
we get the following figure.

\psset{xunit=1cm,yunit=1cm, algebraic=true}
\begin{pspicture}(-1.9,-1.9)(8,4)
\psgrid[gridwidth=0.1pt,gridlabels=0,subgriddiv=2,gridcolor=lightgray]
\psline*[linewidth=2pt,%
linecolor=lightgray]%
{-}(0,0.5)(0.33,0.33)(0,1)(0,0.5)
\psaxes[Dx=1,Dy=1,dx=1,dy=1, subticks=2,tickstyle=buttom]{->}(0,0)(-1.9,-1.9)(8,4)
\rput[l](7.7,-0.3){$z_2$}
\rput[l](-0.4,3.8){$z_3$}
\psplot[linewidth=2pt, linestyle=solid]{0}{1}{-2*x+1}
\psplot[linewidth=2pt, linestyle=solid]{0}{2}{-0.5*x+0.5}
\thinlines
\put(1,-1){\vector(-2,-1){0.5}}
\put(1,0){\vector(1,2){0.3}}
\end{pspicture}

With $ w_1  \leq  2z_2 -1$ and $z_2 \in [0,\frac{1}{3}]$ we get
\[
w_1 \leq  \frac{2}{3} -1 =-\frac{1}{3}<0.
\]
This violates $w_1 \geq 0$. So, Case 4 gives no contribution to the solution set. 

Case 6: $w_1=0$, $w_2=0$, $z_3=0$, $z_1\geq 0$, $z_2\geq 0$, $w_3\geq 0$. Then, the inequalities reduce to
\[
\begin{array}{rclcrcl}
 -5z_1-2z_2  & \leq & -1 & \quad \mbox{and} \quad & -2 & \leq & -4z_1+z_2   \\[1ex]
 -2z_1-3z_2  & \leq & 1 & \quad \mbox{and} \quad & -1 & \leq & z_1-2z_2  \\[1ex]
w_3 -z_2  & \leq & -1 & \quad \mbox{and} \quad & -2 & \leq & w_3 -z_2 ,  
\end{array}
\]
whence
\[
\begin{array}{rclcrcl}
\displaystyle{-\frac{5}{2}z_1+\frac{1}{2}} & \leq & z_2 & \quad \mbox{and} \quad & 4z_1 -2 & \leq & z_2   \\[2ex]
\displaystyle{-\frac{2}{3}z_1-\frac{1}{3}} & \leq & z_2 & \quad \mbox{and} \quad & z_2 & \leq & \displaystyle{\frac{1}{2}z_1+\frac{1}{2}}  \\[2ex]
w_3 +1 & \leq & z_2 & \quad \mbox{and} \quad & z_2 & \leq & w_3 +2.  
\end{array}
\]
Regarding the first four inequalities together with $z_1\geq 0$, $z_2\geq 0$ we get the following polygon

\psset{xunit=4cm,yunit=4cm, algebraic=true}
\begin{pspicture}(-0.4,-0.4)(1.5,1.5)
\psgrid[gridwidth=0.1pt,gridlabels=0,subgriddiv=2,gridcolor=lightgray]
\psaxes[Dx=1,Dy=1,dx=1,dy=1, subticks=2,tickstyle=buttom]{->}(0,0)(-0.4,-0.4)(1.5,1.5)
\rput[l](1.3,-0.1){$z_1$}
\rput[l](-0.1,1.3){$z_2$}
\psplot[linewidth=2pt, linestyle=solid]{0.5}{0.714}{4*x-2}
\psplot[linewidth=2pt, linestyle=solid]{0}{0.714}{0.5*x+0.5}
\psplot[linewidth=2pt, linestyle=solid]{0}{0.2}{-2.5*x+0.5}
\psplot[linewidth=2pt, linestyle=solid]{0.2}{0.5}{0}
\end{pspicture}

with its supremum $(\frac{5}{7}/\frac{6}{7})$. However, $z_2 \leq \frac{6}{7}$ and $ w_3 +1  \leq  z_2 $ is a contradiction. So, Case 6 gives no contribution to the solution set.

Case 7: $w_1=0$, $z_2=0$, $w_3=0$, $z_1\geq 0$, $w_2\geq 0$, $z_3\geq 0$. Then, the inequalities reduce to
\[
\begin{array}{rclcrcl}
-5z_1  & \leq & -1 & \quad \mbox{and} \quad & -2 & \leq &  -4z_1   \\[1ex]
w_2 -2z_1 -z_3 & \leq & 1 & \quad \mbox{and} \quad & -1 & \leq & w_2 +z_1 -z_3  \\[1ex]
-2z_3 & \leq & -1 & \quad \mbox{and} \quad & -2 & \leq &  -2z_3, 
\end{array}
\]
whence
\[
\begin{array}{rclcrcl}
\displaystyle{\frac{1}{5}} & \leq & z_1 & \quad \mbox{and} \quad & z_1 & \leq &   \displaystyle{\frac{1}{2}} \\[2ex]
-2z_1 +w_2   -z_3 & \leq & 1 & \quad \mbox{and} \quad & -1 & \leq & z_1 + w_2 -z_3  \\[1ex]
\frac{1}{2} & \leq & z_3 & \quad \mbox{and} \quad & z_3 & \leq & 1.  
\end{array}
\]
On the one hand, the intersection of the plane $-2z_1 +w_2   -z_3 =1 $ with $z_1=\frac{1}{5}$, $z_1=\frac{1}{2}$, $z_3=\frac{1}{2}$, and $z_3=1$, respectively, gives the vertices
\[
(\frac{1}{5}/\frac{19}{10}/\frac{1}{2}), \quad (\frac{1}{5}/\frac{12}{5}/1), \quad (\frac{1}{2}/\frac{5}{2}/\frac{1}{2}), \quad (\frac{1}{2}/3/1).   
\]
On the other hand, the intersection of the plane $-1= z_1 +w_2   -z_3  $ with $z_1=\frac{1}{5}$, $z_1=\frac{1}{2}$, $z_3=\frac{1}{2}$, and $z_3=1$, respectively, gives the vertices
\[
(\frac{1}{5}/-\frac{7}{10}/\frac{1}{2}), \quad (\frac{1}{5}/-\frac{1}{5}/1), \quad (\frac{1}{2}/-1/\frac{1}{2}), \quad (\frac{1}{2}/-\frac{1}{2}/1).   
\]
So, the contribution to the solution set is the rectangle
\[
\left\{
\left(
\begin{array}{c}
z_1 \\[1ex]
z_2 \\[1ex]
z_3
\end{array}
\right) \in  \mathbb{R}^3\, :  \, z_1 \in [\displaystyle{\frac{1}{5}},\displaystyle{\frac{1}{2}}], \, z_2=0,\,   z_3 \in  [\displaystyle{\frac{1}{2}},1] 
\right\}.
\]
Case 8: $w_1=0$, $w_2=0$, $w_3=0$, $z_1\geq 0$, $z_2\geq 0$, $z_3\geq 0$. Then, the inequalities reduce to
\[
\begin{array}{rclcrcl}
-5z_1-2z_2  & \leq & -1 & \quad \mbox{and} \quad & -2 & \leq &  -4z_1+z_2   \\[1ex]
 -2z_1-3z_2 -z_3 & \leq & 1 & \quad \mbox{and} \quad & -1 & \leq & z_1-2z_2 -z_3  \\[1ex]
 -z_2 -2z_3 & \leq & -1 & \quad \mbox{and} \quad & -2 & \leq &  -z_2 -2z_3.  
\end{array}
\]
Considering the first, the second, and the two last conditions we get the following polyhedron\\

\begin{figure}[h]
\unitlength1.0cm
\begin{picture}(12,10.5)
\includegraphics[scale=0.5]{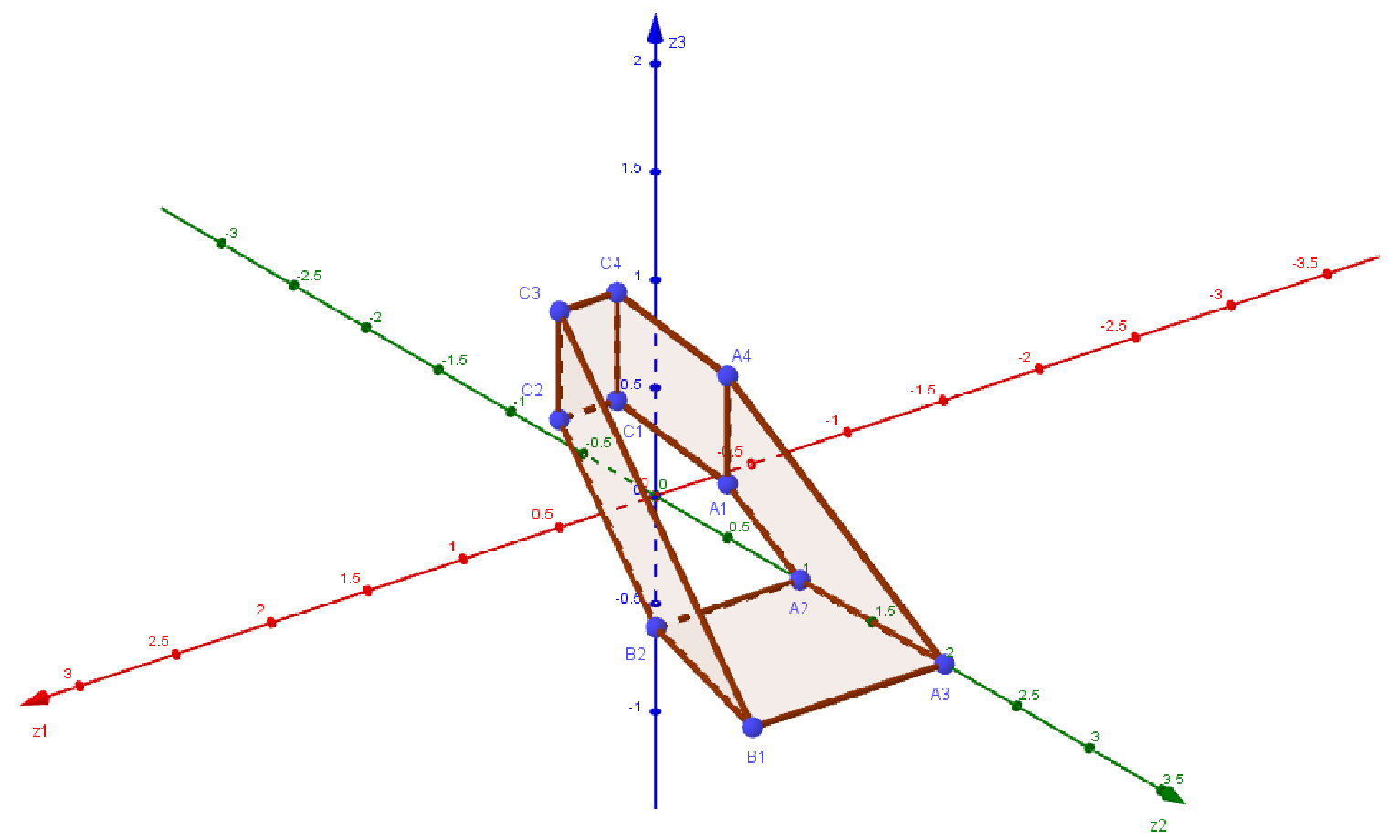}
\end{picture}
\end{figure}

with the vertices
\[
C1(\displaystyle{\frac{1}{5}}/0/\displaystyle{\frac{1}{2}}), \quad
C2(\displaystyle{\frac{1}{2}}/0/\displaystyle{\frac{1}{2}}), \quad
C3(\displaystyle{\frac{1}{2}}/0/1), \quad
C4(\displaystyle{\frac{1}{5}}/0/1),
\]
\[
A1(0/\displaystyle{\frac{1}{2}}/\displaystyle{\frac{1}{4}}), \quad
A2(0/1/0), \quad
A3(0/2/0), \quad
A4(0/\displaystyle{\frac{1}{2}}/\displaystyle{\frac{3}{4}}),
\]
and
\[
B1(1/2/0), \quad
B2(\displaystyle{\frac{3}{4}}/1/0).
\]
The vertices $C1$, $C4$, $A1$, and $A4$ belong to the plane
\[
5z_1+2z_2=1.
\]
The vertices $C2$, $C3$, $B1$, and $B2$ belong to the plane
\[
4z_1-z_2=2.
\]
The vertices $C1$, $C2$, $A1$, $A2$, and $B2$ belong to the plane
\[
z_2+2z_3=1.
\]
The vertices $C3$, $C4$, $A4$, $B1$, and $A3$ belong to the plane
\[
z_2+2z_3=2.
\]
The plane $-2z_1-3z_2-z_3=1$ does not truncate the polyhedron, but the plane
\[
-1  =  z_1-2z_2 -z_3
\]
does. The new vertices are 
\[
 D1(\displaystyle{\frac{7}{10}}/\displaystyle{\frac{4}{5}}/\displaystyle{\frac{1}{10}}), \quad
 D2(\displaystyle{\frac{3}{5}}/\displaystyle{\frac{2}{5}}/\displaystyle{\frac{4}{5}}), \quad
 D3(\displaystyle{\frac{1}{19}}/\displaystyle{\frac{7}{19}}/\displaystyle{\frac{6}{19}}), \quad
 D4(\displaystyle{\frac{3}{19}}/\displaystyle{\frac{2}{19}}/\displaystyle{\frac{18}{19}})
\]
and we get the polyhedron above.\\

\begin{figure}[t]
\unitlength1.0cm
\begin{picture}(12,10.5)
\includegraphics[scale=0.4]{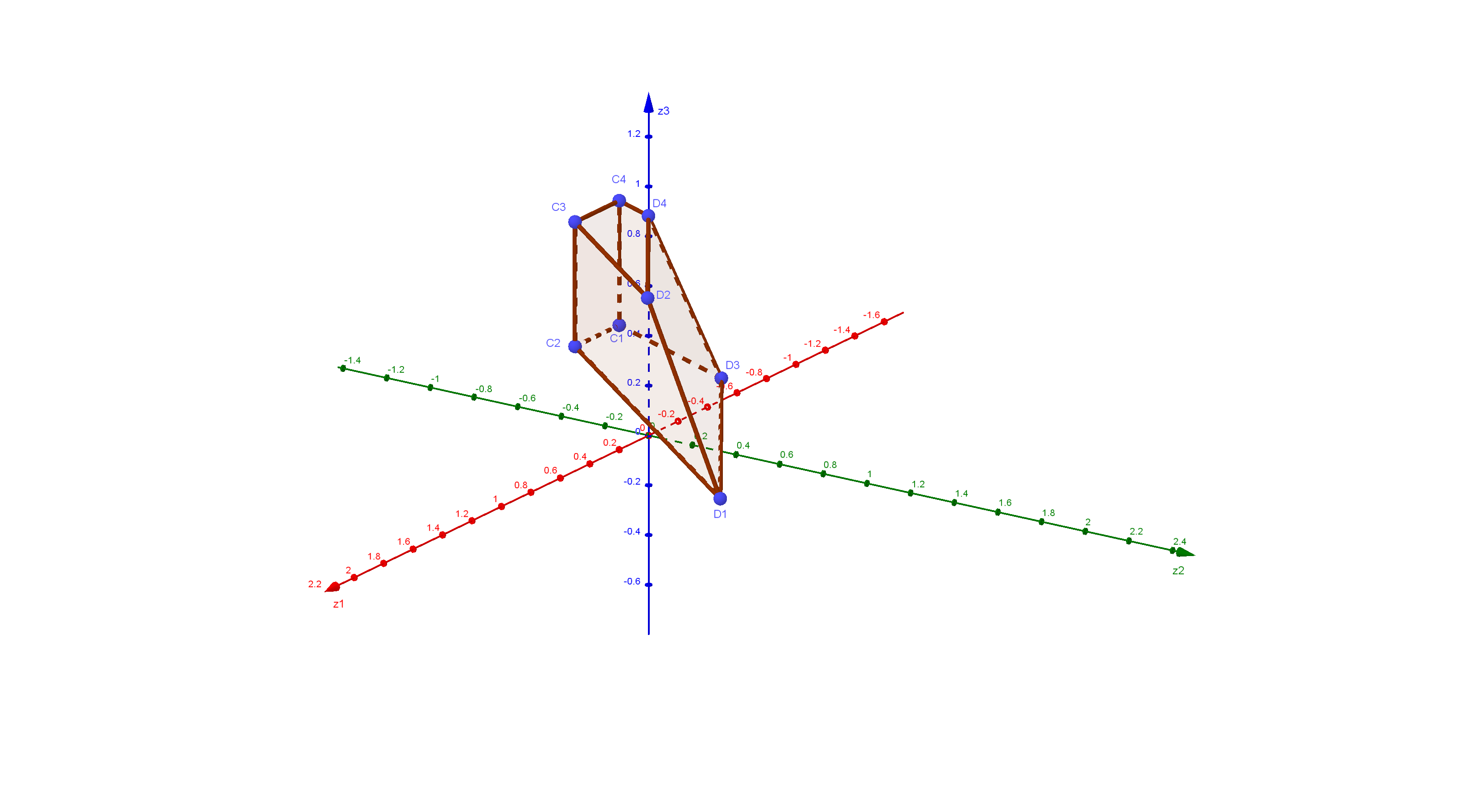}
\end{picture}
\end{figure}

Now we consider the symmetric solution set. First of all we want to note that the rectangle $C1C2C3C4$ must belong to the symmetric solution set. This rectangle was the contribution to the solution set gained by Case 7, where $w_1=0$, $z_2=0$, and $w_3=0$; i.e., for any $z$ corresponding to a point that belongs to the rectangle $C1C2C3C4$ there are $q \in [q]$, $m_{11}
\in [m_{11}]$, $ m_{21} \in [m_{21}]$, and $w_2\geq 0$ such that
\[
q=
\left(
\begin{array}{ccc}
-m_{11} & 0 & 0 \\[1ex]
-m_{21} & 1 & -1 \\[1ex]
0& 0& -2
\end{array}
\right)\cdot
\left(
\begin{array}{c}
z_{1} \\[1ex]
w_{2}  \\[1ex]
z_3
\end{array}
\right).
\]
Choosing $m_{12}=m_{21}$ and choosing $m_{22} \in [m_{22}]$ arbitrarily, then $z \in \Sigma_{sym}([M],[q])$. Next we will show that the vector $z$ corresponding to the 
point $D3(\frac{1}{19}/\frac{7}{19}/\frac{6}{19})$ does not belong to the symmetric solution set. Indeed, we have
\[
-\Big(
[m_{11}] \, \, [m_{12}] \, \, [m_{13}]
\Big) \cdot
\left(
\begin{array}{c}
\frac{1}{19}\\[1ex]
\frac{7}{19}\\[1ex]
\frac{6}{19}
\end{array}
\right)
=-[4,5] \cdot \frac{1}{19} -[-1,2]\cdot  \frac{7}{19} =  [-1,\frac{3}{19}] \cap \underbrace{[-2,-1]}_{=[q_1]} \not= \emptyset.
\]
Necessarily, we get $m_{11}=5$ and $m_{12}=2$. So,
\[
z \in \Sigma_{sym}([M],[q]) \quad \Leftrightarrow \quad m_{21}=2.
\]
However,
\[
q_2 + \Big(
2 \, \, m_{22} \, \, 1
\Big) \cdot
\left(
\begin{array}{c}
\frac{1}{19}\\[1ex]
\frac{7}{19}\\[1ex]
\frac{6}{19}
\end{array}
\right) = q_2 + \frac{7}{19} \cdot m_{22} + \frac{8}{19} \in [-1,1] +\frac{7}{19}\cdot [2,3] + \frac{8}{19}=[\frac{3}{19},\frac{48}{19}] \not\ni 0. 
\]
Therefore, in contrast to Example 6.2.2 $\Sigma_{sym}([M],[q])$ will be smaller than $\Sigma([M],[q])$. The process described in Section 5.2.1 gives
\[
0\leq 3 \cdot z_2^2 + z_2 \cdot z_3 + z_2 -4 \cdot z_1^2 + 2\cdot z_1
\quad
\mbox{and} \quad
0 \leq 5 \cdot z_1^2 -z_1 -2 \cdot z_2^2 -z_2 \cdot z_3 + z_2.
\]
The intersection of the hyperbolic paraboloid
\[
H_2: \quad 0 = 5 \cdot z_1^2 -z_1 -2 \cdot z_2^2 -z_2 \cdot z_3 + z_2
\]
with the quadrangle $C3C4D4D2$ looks like the following figure,
\newline
\includegraphics[scale=1.5]{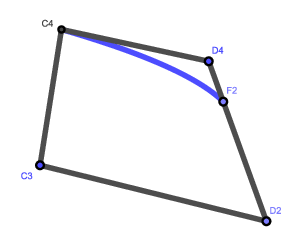}

with
\[
F2(\frac{3}{13}/\frac{2}{13}/ \frac{12}{13}).
\]
The intersection of $H_2$ with the quadrangle $C1C2D1D3$ looks like the following figure,
\newline
\includegraphics[scale=1.3]{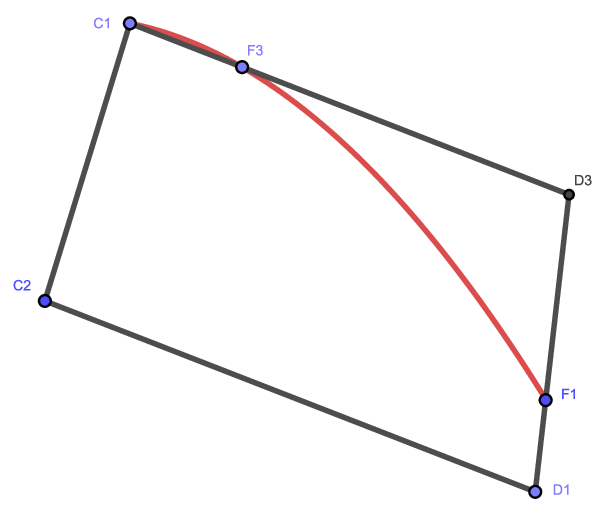}

with
\[
F1(\frac{4}{13}/\frac{7}{13}/ \frac{3}{13}) \quad \mbox{and} \quad
F3(\frac{1}{7}/\frac{1}{7}/ \frac{3}{7}).
\]
The intersection of $H_2$ with the quadrangle $C1C4D4D3$ looks like the following figure,
\newline
\includegraphics[scale=1.2]{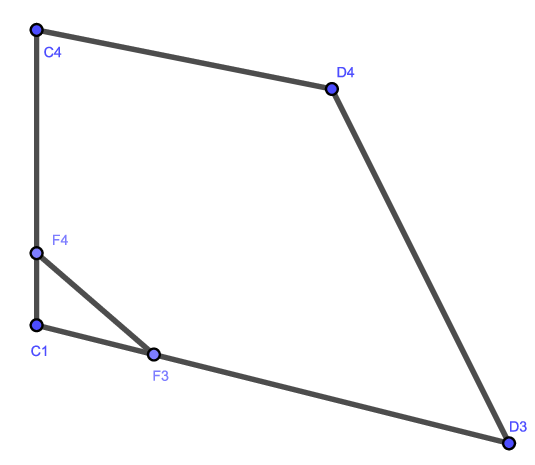}

where the path from $F3$ to $F4$ is a straight line with
\[
F4(\frac{1}{5}/0/ \frac{3}{5}).
\]
Finally, the intersection of $H_2$ with the quadrangle $D1D2D3D4$ looks like the following figure,
\newline
\includegraphics[scale=1.3]{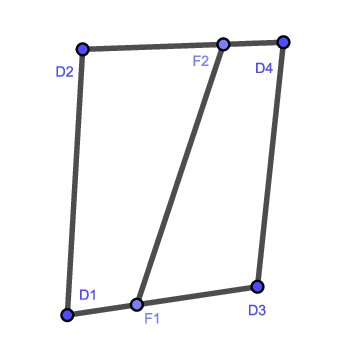}

where the path from $F1$ to $F2$ is a straight line. We do not go into further details and refer to the following link.

\url{https://youtu.be/Umabs6SHK8w}

\end{document}